\documentclass[12pt, a4paper]{amsart}
\usepackage{amscd,amssymb,eucal,eufrak,mathrsfs,amsmath}
\usepackage[hidelinks]{hyperref}
\usepackage[usenames,dvipsnames]{color}

\input xy
\xyoption{all}
\usepackage{epsfig}

\newtheorem{thm}{Theorem}[section]
\newtheorem{cor}[thm]{Corollary}
\newtheorem{lem}[thm]{Lemma}
\newtheorem{prop}[thm]{Proposition}
\newtheorem{fact}[thm]{Fact}

\newtheorem{question}[thm]{Question}

\newtheorem{claim}[thm]{Claim}

\theoremstyle{definition}
\newtheorem{defn}[thm]{Definition}

\theoremstyle{remark}
\newtheorem{remark}[thm]{Remark}

\newtheorem{example}[thm]{Example}

\numberwithin{equation}{section}
\newtheorem{problem}[thm]{Problem}

\newcommand{\delete}[1]{} 
\newcommand{\nt}{\noindent}

\def\eps{{\varepsilon}}
\newcommand{\sk}{\vskip 0.2cm}

\newcommand{\bsk}{\vskip 0.4cm}

\newcommand{\nl}{\newline}
\newcommand{\ben}{\begin{enumerate}}

\newcommand{\een}{\end{enumerate}}
\newcommand{\bit}{\begin{itemize}}

\newcommand{\eit}{\end{itemize}}

\newcommand{\nbd}{neighborhood }
\newcommand{\wrt}{with respect to }
\newcommand{\rest}{\upharpoonright}

\def\R {{\mathbb R}}
\def\N {{\mathbb N}}
\def\Z {{\mathbb Z}}
\def\Q {{\mathbb Q}}

\def\U{{\mathbb U}}

\def\L{{\mathcal{L}}}

\def\norm#1{\left\Vert#1\right\Vert}
\def\RUC{{\hbox{RUC\,}^b}}

\def\Iso{{\mathrm{Iso}}\,}
\def\Aut{{\mathrm Aut}\,}

\newcommand{\eva}{\rm{eva}}

\def\QED{\nobreak\quad\ifmmode\roman{Q.E.D.}\else{\rm Q.E.D.}\fi}

\def\a{\alpha}
\def\om{\omega}
\def\Om{\Omega}
\def\g{\gamma}
\def\SUC{\mathrm{SUC}} 
\def\WAP{\mathrm{WAP}} 
\def\UC{\mathrm{UC}} 
\def\LE{\mathrm{LE}} 
\def\Asp{\mathrm{Asp}}
\def\RUC{\mathrm{RUC}} 
\def\LUC{\mathrm{LUC}} 
\def\SL{\mathrm{SL}} 

\newcommand{\ga}{\gamma}

\newcommand{\Del}{\Delta}
\newcommand{\br}{\vspace{3 mm}}

\newcommand{\cls}{\rm{cl\,}}

\newcommand{\adj}{\rm{adj}}
\newcommand{\Acal}{\mathcal{A}}

\begin{document}

\baselineskip=17pt     

\title
[]{New algebras of functions on topological groups arising from
$G$-spaces}

\author[]{E. Glasner}
\address{Department of Mathematics, Tel-Aviv University,
Ramat Aviv, Israel} \email{glasner@math.tau.ac.il}
\urladdr{http://www.math.tau.ac.il/$^\sim$glasner}

\author[]{M. Megrelishvili}
\address{Department of Mathematics,
Bar-Ilan University, 52900 Ramat-Gan, Israel}
\email{megereli@math.biu.ac.il}
\urladdr{http://www.math.biu.ac.il/$^\sim$megereli}

\date{January 2025} 


\begin{abstract}    
For a topological group $G$ we introduce the algebra $\SUC(G)$ of
\emph{strongly uniformly continuous} functions. 
It contains the algebra $\WAP(G)$ of weakly almost periodic
functions as well as the algebras $\LE(G)$ and $\Asp(G)$ of locally
equicontinuous and Asplund functions respectively. For the Polish
groups of order preserving homeomorphisms of the unit interval and
of isometries of the Urysohn space of diameter 1, 
$\SUC(G)$ is trivial. 
We study the Roelcke algebra (= $\UC(G)$ = right and left
uniformly continuous functions) and SUC compactifications of the
groups $S(\N)$, of permutations of a countable set, and $H(C)$, the
group of homeomorphisms of the Cantor set. For the first group we
show that $\WAP(G)=\SUC(G)=\UC(G)$ and also provide a concrete
description of the corresponding metrizable (in fact Cantor)
semitopological semigroup compactification. For the second group,
in contrast, 
$\SUC(G)$ is properly contained in $\UC(G)$ 
and for this group $\UC(G)$ does not yield a right 
topological semigroup compactification. 

We introduce the notion of fixed point on a
class $P$ of flows ($P$-$fpp$) and study in particular groups 
which are SUC-amenable and groups with the $\SUC$-$fpp$ (SUC-extreme amenability).    
We show that every Polish group $G$ with metrizable $M(G)$ is SUC-amenable and if, in addition, $M(G)$ is proximal, then $G$ is  SUC-extremely amenable. 
\end{abstract}

\thanks{Research partially supported by BSF (Binational USA-Israel)
grant no. 2006119.}

\maketitle
\tableofcontents
\setcounter{tocdepth}{1}

\section{Introduction}
\label{s:intro}

In this paper we introduce the property of {\it Strong Uniform
Continuity} (in short: SUC) of $G$-spaces and the associated
notion of SUC functions. For every compact $G$-space $X$ the
corresponding orbit maps $\widetilde{x} \colon G \to X, \ g \mapsto gx$
are right uniformly continuous for every $x \in X$. If all the
maps $\{\widetilde{x}\}_{x \in X}$ are also left uniformly
continuous then we say that $X$ is SUC. Every right uniformly
continuous bounded real valued function $f \colon G \to \R$ comes from
some compact $G$-space $X$. That is, there exist a compact
$G$-space $X$, a continuous function $F \colon X \to \R$, and a point $x_0 \in X$ such that $f=F \circ \widetilde{x_0}$. We say that $f$ is SUC if it comes from a compact $G$-space $X$ which is SUC. Denote by $\SUC(G)$ the corresponding class of functions on $G$. The class $\SUC(G)$ forms a uniformly closed $G$-invariant
subalgebra of the algebra $\UC(G):=\RUC(G) \cap \LUC(G)$ of (right
and left) uniformly continuous functions. Of course we have
$\SUC(G) = \UC(G) = \RUC(G) = \LUC(G)$ 
when $G$ is either discrete or
abelian so that the notion of strong uniform continuity can be
useful only when one deals with non-abelian non-discrete
topological groups. Mostly we will be interested in Polish
non-locally compact large groups, but some of the questions we
study are of interest in the locally compact case as well.

In our recent work \cite{GM} we investigated, among other topics,
the algebras of {\it locally equicontinuous} $\LE(G)$, and {\it
Asplund} functions $\Asp(G)$ on a topological group $G$. The
inclusions 
$\UC(G) \supset \SUC(G) \supset \LE(G) \supset \Asp(G)
\supset \WAP(G)$ hold for an arbitrary topological group $G$. In
the present article we provide a characterization of the elements
of $\SUC(G)$ and $\LE(G)$ in terms of  matrix coefficients for
appropriate Banach representations of $G$ by linear isometries.

Intuitively the dynamical complexity of a function $f \in \RUC(G)$
can be estimated by the \emph{topological complexity} of the
\emph{cyclic} $G$-flow $X_f$ (the pointwise closure of the 
left $G$-orbit $\{gf\}_{g \in G}$ of $f$) treating it as a subset of the Banach space $RUC(G)$. This 
leads to a natural dynamical hierarchy (see Theorem
\ref{r:hierarchy}) where $\SUC(G)$ plays a basic role. In some
sense $\SUC(G)$ is the largest ``nice subalgebra" of $\UC(G)$. It
turns out that $f \in \SUC(G)$ iff $X_f$ is a subset of $\UC(G)$.
Moreover the algebra $\SUC(G)$ is \emph{point-universal} in the
sense of \cite{GM} and every other point-universal subalgebra of
$\UC(G)$ is contained in $\SUC(G)$. Recall that a $G$-algebra
${\mathcal A} \subset \RUC(G)$ is point-universal if and only if
the associated $G$-compactification $G \to G^{{\mathcal A}}$ is a
right topological semigroup compactification.

As an application we conclude that the algebra $\UC(G)$ is
point-universal if and only if it coincides with the algebra
$\SUC(G)$ and that the corresponding {\it Roelcke compactification}
$G \to G^{\UC}$ is in general not a right topological semigroup
compactification of $G$ (in contrast to the compactification $G
\to G^{\SUC}$ determined by the algebra $\SUC(G)$).

For locally compact groups, $\SUC(G)$ contains the subalgebra 
$C_0(G)$ consisting of the functions which vanish at infinity,
and therefore determines the topology of $G$. The
structure of $\SUC(G)$ --- in contrast to $\RUC(G)$ which is always
huge for non-precompact groups --- is ``computable" for several
large groups like: $H_+[0,1]$, $\Iso(\U_1)$ (the isometry group of
the Urysohn space of diameter one $\U_1$), $U(H)$ (the unitary
group on an infinite dimensional Hilbert space), 
$S_{\infty}=S(\N)$ (the Polish infinite symmetric group) and any 
noncompact connected simple Lie group with finite center
(e.g., $\SL_n(\R)$). 
For instance, $\SUC(G)=\WAP(G)$ for $U(H),$ \ 
$S_{\infty}$ and $\SL_n(\R)$. In the first case we use a result of Uspenskij \cite{Us-un} which identifies the Roelcke completion of $U(H)$ as the compact semigroup of contracting operators on the Hilbert space $H$. 
For $S_\infty$ see Section \ref{Sec-S}, and for $\SL_n(\R)$ this 
follows from an old result of Veech \cite{V}.

The group $H_+[0,1]$ of orientation preserving homeomorphisms of
the closed unit interval, endowed with the compact open topology
is a good test case in the class of ``large" yet ``computable"
topological groups. See Section \ref{s:prop} for more details on
this group. In particular recall the result from \cite{Merup}
which shows that $H_+[0,1]$ is WAP-\emph{trivial}: Every weakly
almost periodic function on $G:=H_+[0,1]$ is a constant.
Equivalently, $G$ is \emph{reflexively trivial}, that is, every
continuous representation $G \to \Iso(V)$ where $V$ is a reflexive
Banach space is trivial.

Here we show that $G$ is even ``SUC-trivial"
--- that is, the algebra $\SUC(G)$
(and hence, also the algebras $\LE(G)$ and $\Asp(G)$) consists only
of constant functions --- and that every continuous representation
of $G$ into the group of linear isometries $\Iso(V)$ of an Asplund
Banach space $V$ is trivial. Since in general $\WAP(G) \subset
Asp(G)$ and since every reflexive Banach space is Asplund these
results strengthen the main results of \cite{Merup}.
SUC-triviality implies that every \emph{adjoint continuous} (see
Section \ref{s:repr}) representation is trivial for $H_+[0,1]$. 
The latter fact follows also from a recent unpublished result of
Uspenskij (private communication).


From the WAP-triviality (equivalently, reflexive triviality) 
of $H_+[0,1]$ and results of Uspenskij about $\Iso(\U_1)$, Pestov
deduces in \cite[Corollary 1.4]{Pe-new} the fact that the group
$\Iso(\U_1)$ is also WAP-trivial. Using a similar idea and the
matrix coefficient characterization of SUC one can conclude that
$\Iso(\U_1)$ is SUC-trivial. It is an open question whether the
group $H([0,1]^{\omega})$ is SUC-trivial (or, WAP-trivial).

The above mentioned description of $\SUC(G)$ and $\LE(G)$ in terms
of matrix coefficients (Section \ref{s:mat-SUC}) is nontrivial.
The proof is based on a dynamical modification of a well known
interpolation technique of Davis, Figiel, Johnson and Pelczy\'nski
\cite{DFJP}.

In the last two sections we
study the Roelcke and SUC compactifications of the groups
$S_\infty$ and $H(C)$. For the first group we show that
$\WAP(S_\infty)=\SUC(S_\infty)=\UC(S_\infty)$ and also provide a
concrete description of the corresponding metrizable (in fact
Cantor) semitopological semigroup compactification. 
For the latter group $G:=H(C)$, in contrast, we have
$\SUC(G)\subsetneqq \UC(G)$ from which fact we deduce that the
corresponding Roelcke compactification $G \to G^{\UC}$ is not a
right topological semigroup compactification of $G$.

Finally let us note that although in this work we consider,
for convenience, algebras of real-valued functions,
it seems that there should be no difficulty in extending our
definitions and results to the complex case.

In section \ref{s:amenable} we introduce a notion of 
(amenability) 
extreme amenability with respect to a class of flows. In particular we
examine \emph{extreme SUC-amenability} and \emph{extreme
	SUC-amenable} groups. Namely those groups which have a fixed point
property on compact SUC $G$-spaces. Several natural groups, like
$\mathrm{SL_2}(\R)$, $S_{\infty}$, $H(C)$ (the homeomorphisms group of the 
Cantor set), and $H_+(\mathbb{T})$, which fail to be extremely amenable,
are however extremely SUC-amenable. 

By Theorem \ref{t:GeneralCor} 
if $G$ is a Polish topological group such that the universal minimal $G$-flow $M(G)$ is metrizable and proximal, then $M(G)$ is SUC-trivial (hence, $G$ is extremely SUC-amenable). 
The same result can be derived
by using Theorem \ref{t:stronger}. Furthermore, 
Theorem \ref{t:GeneralCor} 
leads to    
Proposition \ref{t:SUC-amenable} which asserts that every Polish $G$ with metrizable $M(G)$ is SUC-amenable.

\section{Actions and $G$-compactifications}
\label{s:actions}

Unless explicitly stated otherwise, all spaces in this paper are
at least Tychonoff. A (\emph{left}) \emph{action} of a topological
group $G$ on a topological space $X$ 
is defined by a function $\pi \colon G \times X \to X, \ \pi(g,x):=gx$
such that always $g_1(g_2x)=(g_1g_2)x$ and $ex=x$ hold, where
$e=e_G$ is the neutral element of $G$.
Every $x \in X$ defines an \emph{orbit map}
${\tilde x}\colon G \to X, \ g \mapsto gx$.
Also every $g \in G$ induces a \emph{$g$-translation}
$\pi^g\colon X \to X, \ x \mapsto gx$. If the action $\pi$ is
continuous then we say that $X$ is a $G$-\emph{space} (or a
$G$-system or a $G$-\emph{flow}).
Sometimes we denote it as a pair $(G,X)$. If the orbit
$Gx_0$ of $x_0$ is dense in $X$ for some $x_0 \in X$ then
the $G$-space $X$ is \emph{point transitive}
(or just \emph{transitive}) and the point $x_0$ is a \emph{transitive point}. 
If $X$ in addition is compact then the pair $(X,x_0)$ is said to
be a \emph{pointed system} or a {\it $G$-ambit}. If every point
$x$ in a compact $G$-space $X$ is transitive then $X$ is said to
be \emph{minimal}.

Let $G$ act on $X_1$ and on $X_2$. A continuous map
$f\colon X_1 \to X_2$ is a \emph{$G$-map}
(or a {\em homomorphism} of dynamical systems)
if $f(gx)=gf(x)$ for every $(g,x) \in G \times X_1$.

A {\it right action} $X \times G \to X$ can be defined analogously.
If $G^{op}$ is the {\it opposite group} of $G$ with the same
topology then the right $G$-space $(X,G)$ can be treated as a left
$G^{op}$-space $(G^{op},X)$ (and vice versa). A map $h\colon G_1 \to
G_2$ between two groups is a {\it co-homomorphism} (or, an
\emph{anti-homomorphism}) if $h(g_1g_2)=h(g_2)h(g_1)$. This
happens iff $h\colon G_1^{op} \to G_2$ (the same assignment) is a
homomorphism.

The Banach algebra (under the supremum norm) of all continuous
real valued  bounded functions on a topological space $X$ will be
denoted by $C(X)$.
 Let $(G,X)$ be a left (not necessarily compact) $G$-space.
 Then it induces the right
action $C(X)\times G \to C(X)$, with $(fg)(x)=f(gx)$, and the
corresponding co-homomorphism $h\colon G \to \Iso(C(X))$. While the
$g$-translations $C(X) \to C(X)$ (being isometric) are continuous,
the orbit maps ${\tilde f}\colon G \to C(X), \ g \mapsto fg$ are not
necessarily continuous. The function $f \in C(X)$ is {\it right
uniformly continuous} if the orbit map $G \to C(X), \ g \mapsto
fg$ is norm continuous. The set $\RUC(X)$ of all right uniformly
continuous functions on $X$ is a uniformly closed $G$-invariant
subalgebra of $C(X)$. Here and in the sequel ``subalgebra" means a
uniformly closed unital (containing the constants) subalgebra. A
``$G$-subalgebra" is an algebra which is invariant under the
natural right action of $G$.

Every topological group $G$ can be treated as a $G$-space under
the left regular action of $G$ on itself. In this particular case
$f \in \RUC(G)$ iff $f$ is uniformly continuous with respect to the
{\it right uniform structure} $\mathcal {R}$ on $G$
(furthermore this is also true for coset $G$-spaces $G/H$).

Thus, $f\in \RUC(G)$ iff for every $\eps>0$ there exists a neighborhood 
$V$ of
the identity element $e\in G$ such that $\sup_{g\in
G}|f(vg)-f(g)|<\eps$ for every $v \in V$.

Analogously one defines {\it right translations} $(gf)(x):=f(xg)$,
and the algebra $\LUC(G)$ of {\it left uniformly continuous}
functions. These are the functions which are uniformly continuous
with respect to the {\it left uniform structure} $\mathcal {L}$ on
$G$.

A $G$-{\it compactification} of a $G$-space $X$ is a
$G$-map $\nu\colon X \to Y$ into a compact $G$-space $Y$ with
${\cls}\nu(X)=Y$. A compactification is {\it proper} when $\nu$ is
a topological embedding. Given a compact $G$-space $X$ and a point
$x_0\in X$ the map $\nu\colon G \to X$ defined by $\nu(g)=gx_0$ is a
compactification of the $G$-space $G$ (the left regular action) in
the orbit closure
${\cls}Gx_0\subset X$.

We say that a $G$-compactification $\nu\colon G \to S$ of $X:=G$
is a {\it right topological semigroup
compactification} of $G$ if $S$ is a {\it right topological
semigroup} (that is, $S$ is a compact semigroup such that for
every $p \in S$ the map $S \to S, s \mapsto sp$ is continuous) and
$\nu$ is a homomorphism of semigroups.

There exists a canonical 1-1 correspondence between the
$G$-compactifications of $X$ and $G$-subalgebras of $\RUC(X)$ (see
for example \cite{Vr-oldpaper}). 
The compactification $\nu\colon X \to
Y$ induces an isometric $G$-embedding of $G$-algebras
$
j_{\nu} \colon C(Y)=\RUC(Y) \hookrightarrow \RUC(X), \ \ \phi
\mapsto \phi \circ \nu
$
and the algebra $\Acal_{\nu}$ is defined as the image
$j_{\nu}(C(Y))$. Conversely, if $\Acal$ is a $G$-subalgebra of
$\RUC(X)$, then denote by $X^{\Acal}$ or by $|\Acal|$ the
corresponding Gelfand space treating it as a weak star compact
subset of the dual space $\Acal^*$. It has a structure of a
$G$-space $(G,|\Acal|)$ and the natural map $\nu_{\Acal}\colon X \to
X^{\Acal}, \hskip 0.3cm x \mapsto eva_x$, where
$eva_x(\varphi):=\varphi(x)$, is the evaluation at $x$ (a
multiplicative functional), defines a $G$-compactification. If
$\nu_1\colon X \to X^{\Acal_1}$ and $\nu_2\colon X \to X^{\Acal_2}$ are two
$G$-compactifications then $\Acal_{\nu_1} \subset \Acal_{\nu_2}$
iff $\nu_1= \a \circ \nu_2$ for some
$G$-map $\a\colon X^{\Acal_2} \to X^{\Acal_1}$. The algebra
$\Acal$ determines the compactification $\nu_{\Acal}$ uniquely, up
to the equivalence of $G$-compactifications. The $G$-algebra
$\RUC(X)$ defines the corresponding Gelfand space $|\RUC(X)|$, which
we denote by $\beta_G X$, and the {\it maximal
$G$-compactification} $i_{\beta}\colon X \to \beta_G X$. Note that this
map may not be an embedding even for Polish $G$ and $X$ (see
\cite{Me-ex}); it follows that there is no proper
$G$-compactification for such $X$. If $X$ is a compact $G$-space
then $\beta_G X$ can be identified with $X$ and $C(X)=\RUC(X)$.

Denote by $G^{\RUC}$ the Gelfand space of the $G$-algebra $\RUC(G)$.
The canonical embedding $u\colon G \to G^{\RUC}$ defines the {\it
greatest ambit} $(G^{\RUC}, u(e))$ of $G$.

It is easy to see that the intersection $\UC(G):=\RUC(G) \cap
\LUC(G)$ is a left and right $G$-invariant closed subalgebra of
$\RUC(G)$. We denote the corresponding compactification by
$G^{\UC}$. Denote by $\mathcal {L} \bigwedge \mathcal {R}$ the
\emph{lower uniformity} of $G$. It is the infimum (greatest lower
bound) of the left and right uniformities on $G$;
we call it the {\it Roelcke uniformity}. Clearly, for every
bounded function $f \colon G \to \R$ we have $f \in \UC(G)$ iff $f\colon
(G,\mathcal {L} \bigwedge \mathcal {R}) \to \R$ is uniformly
continuous. Recall the following important
fact (in general the infimum $\mu_1 \bigwedge \mu_2$ of two compatible
uniform structures on a topological space $X$ need not be
compatible with the topology of $X$).

\begin{lem} \label{l:Roelcke} \ 
\ben
\item \emph{(Roelcke-Dierolf \cite{RD})}
For every topological group $G$ the Roelcke uniform structure
$\mathcal {L} \bigwedge \mathcal {R}$ generates the given topology
of $G$.
\item
For every topological group $G$ the algebra $UC(G)$
separates points from closed subsets in $G$. \een
\end{lem}
\begin{proof} (1)  
See Roelcke-Dierolf \cite[Proposition 2.5]{RD}.

(2) Follows from (1).
\end{proof}

By a {\em uniform\/} $G$-space $(X,\mu)$ we mean a $G$-space
$(X,\tau)$ where $\tau$ is the (completely regular) topology
defined by the uniform structure $\mu$
and the $g$-translations ($g \in G$) are uniform isomorphisms.

Let $X:=(X, \mu)$ be a uniform $G$-space. A point $x_0\in
X$ is a {\em point of equicontinuity\/} (notation: $x_0 \in Eq_X$)
if for every entourage $\eps\in \mu$, there is a neighborhood $O$
of $x_0$ such that $(g x_0,g x)\in \eps$ for every $x\in O$ and
$g\in G$. The $G$-space $X$ is {\it equicontinuous} if $Eq_X=X$.
$(X,\mu)$ is {\it uniformly equicontinuous} if for every $\eps \in
\mu$ there is $\delta \in \mu$ such that $(gx,gy) \in \eps$ for
every $g \in G$ and $(x,y) \in \delta$. For compact $X$ (equipped
with the unique compatible uniformity), equicontinuity and uniform
equicontinuity coincide. Compact (uniformly) equicontinuous
$G$-space $X$ is also said to be \emph{Almost Periodic} (in short:
AP); see also Section \ref{s:hierarchy}. If $Eq_X$ is dense in $X$
then $(X,\mu)$ is said to be an {\it almost equicontinuous} (AE)
$G$-space \cite{AAB2}.

\sk

The following definition is standard (for more details see for
example \cite{GM}).

\begin{defn} \label{d:coming} \ 
\ben
\item
A function $f\in C(X)$ on a $G$-space $X$ {\em comes from\/} a
compact $G$-system $Y$ if there exist a $G$-compactification $\nu:
X \to Y$ (so, $\nu$ is onto if $X$ is compact) and a function
$F\in C(Y)$ such that $f=F \circ \nu$ (i.e., if $f \in
\Acal_{\nu}$).
Then necessarily, $f \in \RUC(X)$.
\begin{equation*}
\xymatrix { X \ar[dr]_{f} \ar[r]^{\nu} & Y
\ar[d]^{F} \\
  & \R }
\end{equation*}
\item
A function $f\in \RUC(G)$ {\em comes from\/} a pointed system
$(Y,y_0)$ if for some continuous function $F \in C(Y)$ we have
$f(g)=F(g y_0),\ \forall g\in G$. Notation:  $f \in \Acal(Y,y_0)$.
Defining $\nu: X=G \to Y$ by $\nu(g)=gy_0$ observe that this is
indeed a particular case of \ref{d:coming}.1.
\item
Let $\Gamma$ be a class of compact $G$-spaces. For a $G$-space $X$
denote by $\Gamma(X)$ the class of all functions on $X$ which come
from a $G$-compactification $\nu\colon X \to Y$ where the $G$-system
$Y$ belongs to $\Gamma$.
 \een
\end{defn}

Let P be a class of compact $G$-spaces which is preserved by
$G$-isomorphisms, products and closed $G$-subspaces.
It is well known (see for example, \cite[Proposition 2.9]{GM})
that for every $G$-space $X$ there exists a \emph{universal}
(maximal)
$G$-compac{\allowbreak}tification $X \to X^{{\mathcal P}}$ such that
$X^{{\mathcal P}}$ lies in P. More precisely, for every (not
necessarily compact) $G$-space $X$ denote by ${\mathcal P} \subset
C(X)$ the collection of functions coming from $G$-spaces having
property $P$. Then ${\mathcal P}$ is a uniformly closed,
$G$-invariant subalgebra of $\RUC(X)$ and the maximal
$G$-compactification of $X$ with property $P$ is the corresponding
Gelfand space $X^{{\mathcal P}}:=|{\mathcal P}|$.
If $X$ is compact then $(G,X^{{\mathcal P}})$ is the maximum
factor of $(G,X)$ with property P. In particular let P be one of
the following natural classes of compact $G$-spaces: a) Almost
Periodic (=equicontinuous); b) Weakly Almost Periodic; c)
Hereditarily Not Sensitive; d) Locally Equicontinuous; e) all
compact $G$-spaces. Then,
in this way the following maximal (in
the corresponding class) $G$-compactifications are:
a) $G^{AP}$; b) $G^{WAP}$; c) $G^{Asp}$; d) $G^{LE}$; e)
$G^{\RUC}$.

For undefined concepts and more details see Section
\ref{s:hierarchy} and also \cite{GM}.

\section{Cyclic $G$-systems and point-universality}
\label{s:cyclic}

Here we give some background material about \emph{cyclic} compact
$G$-systems $X_f$ defined for $f \in \RUC(X)$. These $G$-spaces
play a significant role in many aspects of topological dynamics
and are well known at least for the particular case of $X:=G$. We
mostly use the presentation and results of \cite{GM} (see also
\cite{Vr}).

As a first motivation note a simple fact about Definition
\ref{d:coming}. For every $G$-space $X$ a function $f\colon X \to \R$
lies in $\RUC(X)$ iff it comes from a compact $G$-flow $Y$. We can
choose $Y$ via the maximal $G$-compactification $G \to \beta_G X=Y$
of $X$. This is the largest possibility in this setting. Among all
possible $G$-compactifications $\nu\colon X \to Y$ of $X$ such that $f$
comes from $(\nu,Y)$ there exists also the smallest one. Take
simply the smallest $G$-subalgebra ${\mathcal A_f}$ of $\RUC(X)$
generated by the orbit $fG$ of $f$ in $\RUC(X)$. Denote by $X_f$
the Gelfand space $|\mathcal A_f|=X^{\mathcal A_f}$ of the algebra
${\mathcal A_f}$. Then the corresponding $G$-compactification $X
\to Y:=X_f$ is the desired one. We call ${\mathcal A_f}$ and $X_f$
the \emph{cyclic $G$-algebra} and \emph{cyclic $G$-system} of $f$,
respectively. Next we provide an alternative construction and some
basic properties of $X_f$.

Let $X$ be a (not necessarily compact) $G$-space. Given $f \in
\RUC(X)$ let $I=[-\|f\|,\|f\|] \subset \R$ and $\Om=I^G$, the
product space equipped with the compact product topology. We let
$G$ act on $\Om$ by $g\om(h)=\om(hg),\ g,h\in G$.
 Define the continuous map
$$
f_{\sharp}\colon X \to \Om, \hskip 0.3cm f_{\sharp}(x)(g)=f(gx)
$$
and the closure $X_f:= {\cls} (f_{\sharp}(X))$ in $\Om$. Note that
$X_f = f_{\sharp}(X)$ whenever $X$ is compact.

Denoting the unique continuous extension of $f$ to $\beta_G X$ by
$\tilde f$ (it exists because $f \in \RUC(X)$) we now define a map
$$
\psi\colon \beta_G X \to X_f \quad \text{by} \quad \psi(y)(g)= \tilde
f(gy), \qquad y\in \beta_G X, \ g\in G.
$$
Let $pr_e\colon \Om \to \R$ denote the  projection of $\Om=I^G$ onto
the $e$-coordinate and let $F_e:=pr_e \rest_{X_f}\colon X_f \to \R$ be
its restriction to $X_f$. Thus, $F_e(\om):=\om(e)$ for every $\om
\in X_f$.

As before denote by $\Acal_f$ the smallest (closed and unital, of
course) $G$-invariant subalgebra of $\RUC(X)$ which contains $f$.
There is then a naturally defined $G$-action on the Gelfand space
$X^{\Acal_f}=|\Acal_f|$ and a $G$-compactification (morphism of
dynamical systems if $X$ is compact) $\pi_f\colon X \to |\Acal_f|$.
Next consider the map $\pi\colon \beta_G X \to |\Acal_f|$, the
canonical extension of $\pi_f$ induced by the inclusion $\Acal_f
\subset \RUC(X)$.

The action of $G$ on $\Omega$ is not in general continuous.
However, the restricted action on $X_f$ is continuous for every $f
\in \RUC(X)$. This follows from the second assertion of the next
fact.

\begin{prop} \label{f:X_f}
\emph{(See for example \cite{GM})}
\begin{enumerate}
\item
Each $\om \in X_f$ is an element of $\RUC(G)$. That is, $X_f
\subset \RUC(G)$.
\item
The map $\psi\colon \beta_G X \to X_f$ is a continuous homomorphism of
$G$-systems. The dynamical system $(G,|\Acal_f|)$ is isomorphic to
$(G,X_f)$ and the diagram
\begin{equation*}
\xymatrix { X \ar[d]_{\pi_f} \ar[r]^{i_{\beta}} \ar[dr]^-{\pi}  &
\beta_G X \ar[dl]^{f_{\sharp}} \ar[d]^{\psi} \ar[dr]^{\tilde{f}} & \\
|\Acal_f| \ar[r]  &  X_f \ar[l] \ar[r]^{F_e} & \R & }
\end{equation*}
commutes.
\item $f=F_e \circ f_{\sharp}$.
Thus every $f \in \RUC(X)$
comes from the system $X_f$. Moreover, if $f$
comes from a system $Y$ and a $G$-compactification $\nu\colon X \to Y$
then there exists a homomorphism $\a\colon Y \to
X_f$ such that $f_{\sharp}=\a \circ \nu$.
In particular, $f \in \Acal_f \subset \Acal_{\nu}$.
\end{enumerate}
\end{prop}

If $X:=G$ with the usual left action then $X_f$ is the pointwise
closure of the $G$-orbit $Gf:=\{gf\}_{g \in G}$ of $f$ in
$\RUC(G)$. Hence $(X_f,f)$ is a transitive pointed $G$-system.

As expected by the construction the cyclic $G$-systems $X_f$
provide ``building blocks" for compact $G$-spaces. That is, every
compact $G$-space can be embedded into the $G$-product of
$G$-spaces $X_f$.

Let us say that a topological group $G$ is \emph{uniformly
Lindel\"{o}f}
if for every nonempty open subset $O \subset G$ countably many
translates $g_nO$ cover $G$ (there are several alternative names
for this notion: $\om$-bounded, $\omega$-bounded, $\om$-narrow,
$\om$-precompact). It is well known that $G$ is uniformly
Lindel\"{o}f iff $G$ is a topological subgroup in a product of
second countable groups. When $G$ is uniformly Lindel\"{o}f (e.g.
when $G$ is second countable) the compactum $X_f$ is metrizable.

The question ``when is $X_f$ a subset of $UC(G)$?" provides
another motivation for introducing the notion of strongly
uniformly continuous (SUC) functions (see Definition \ref{d:SUC}
and Theorem \ref{p:X_fin}).

The {\em enveloping (or Ellis) semigroup\/} $E=E(X)$ of a compact
$G$-space $X$ is defined as the closure in $X^X$ (with its compact
pointwise convergence topology) of the set $\{\pi^g\colon X \to X\}_{g
\in G}$ of translations considered as a subset of $X^X$. With the
operation of composition of maps this is a right topological
semigroup. Moreover, the map
$$j=j_X\colon G \to E(X), \ g \mapsto \pi^g$$ is a right topological
semigroup compactification of $G$. The compact space $E(X)$
becomes a $G$-space with respect to the natural action
$$G \times E(X) \to E(X), \ \ \ \ (gp)(x)=gp(x).$$

Moreover the pointed $G$-system $(E(X), j(e))$ is point-universal
in the following sense.

\begin{defn} (\cite{GM})
A pointed $G$-system $(X,x_0)$ is {\em point-universal\/} if it
has the property that for every $x\in X$ there is a homomorphism
$\pi_x \colon (X,x_0)\to ({\cls}(Gx),x)$. The $G$-subalgebra $\Acal
\subset \RUC(G)$ is said to be {\em point-universal\/} if the
corresponding $G$-ambit $(G^{\Acal}, u_{\Acal}(e))$ is
point-universal.
\end{defn}

We will use the following characterization of point-universality
from \cite{GM}.

\begin{lem} \label{l:pu}
Let $(X,x_0)$ be a transitive compact $G$-system. The following
conditions are equivalent: \ben
\item
The system $(X,x_0)$ is point-universal.
\item
The orbit map $G \to X, \ g \mapsto gx_0$ is a right topological
semigroup compactification of $G$.
\item
$(X,x_0)$ is $G$-isomorphic to its enveloping semigroup
$(E(X),j(e))$.
\item
$\Acal(X,x_0) = \bigcup_{x \in X} \Acal (
{\cls}(Gx), x)$.
\item
$X_f \subset \Acal(X,x_0)$ for every $f \in \Acal(X,x_0)$ (where
$\Acal(X,x_0)$ is the corresponding subalgebra of $\RUC(G)$ coming
from the $G$-compactification $\nu\colon G \to X, \ \nu(g)=gx_0$). \een
\end{lem}

In particular, for every right topological semigroup
compactification $\nu\colon G \to S$ the pointed $G$-space $(S,
\nu(e))$ is point-universal.
For other properties of point-universality see \cite{GM} and
Remark \ref{r:introverted}.

\section{Strong Uniform Continuity}
\label{s:SUC}

Let $G$ be a topological group. As before denote by $\mathcal{L}$
and $\mathcal{R}$ the left and right uniformities on $G$. We start
with a simple observation.

\begin{lem} \label{l:right}
For every compact $G$-space $X$ the corresponding orbit maps
$$
\widetilde{x}\colon (G, \mathcal{R}) \to (X,\mu_X), \ \ \ g \mapsto gx
$$
are uniformly
continuous for every $x \in X$.
\end{lem}
\begin{proof}
Let $V$ be an open neighborhood of the diagonal $\Del\subset X
\times X$. In order to obtain a contradiction suppose that for
every neighborhood $U$ of $e\in G$ there is $u_{U}\in U$ and
$x_{U}\in X$ such that $(x_{U},u_{U}x_{U}) \not\in V$. For a
convergent subnet we have $\lim (x_{U},u_{U}x_{U}) = (x,x')\not\in
V$ contradicting the joint continuity of the $G$-action.
\end{proof}

In general, for non-commutative groups, one cannot replace
$\mathcal{R}$ by the left uniformity $\mathcal{L}$ (see Remark
\ref{r:abelian}). This leads us to the following definition.

\begin{defn} \label{d:SUC}
Let $G$ be a topological group.
\begin{enumerate}
\item
We say that a uniform $G$-space $(X,\mu)$ is {\it strongly
uniformly continuous at} $x_0 \in X$ (notation: $x_0 \in \SUC_X$)
if the orbit map $\widetilde{x_0}\colon G \to X, \ g \mapsto gx_0$ is
$(\mathcal{L}, \mu)$-uniformly continuous. Precisely, this means
that for every $\varepsilon \in \mu$ there exists a neighborhood
$U$ of $e \in G$ such that
$$(gux_0,gx_0) \in  \eps$$ for every $g\in G$ and every $u \in U$.
If $\SUC_X=X$ we say that $X$ is {\it strongly uniformly
continuous}. 
\item
If $X$ is a compact $G$-space then there exists a unique
compatible uniformity $\mu_X$ on $X$. So, SUC$_X$ is well defined.
By Lemma \ref{l:right} it follows that a compact $G$-space $X$ is
SUC at $x_0$ iff $\widetilde{x_0}\colon G \to X$ is $(\mathcal {L}
\bigwedge \mathcal {R}, \mu_X)$-uniformly continuous. We let SUC
denote the class of all compact $G$-systems such that $(X,\mu_X)$
is SUC.
\item A function $f \in C(X)$ is {\it strongly
uniformly continuous} (notation: $f \in \SUC(X)$) if it comes from
a SUC compact dynamical system.
\item
Let $x_1$ and $x_2$ be points of a $G$-space $X$. Write $x_1
\stackrel{\SUC} {\sim} x_2$ if these points cannot be separated by
a SUC function on $X$. Equivalently, this means that these points
have the same images under the \emph{universal SUC
compactification} $G$-map $X \to X^{\SUC}$ (see Lemma \ref{p:qLE}.1
below).
\end{enumerate}
\end{defn}

\begin{lem} \label{l:SUC-s-UC}
$\SUC(G) \subset \UC(G)$.
\end{lem}
\begin{proof}
Let $f\colon G \to \R$ belong to $\SUC(G)$. Then it comes from a
function $F\colon X \to \R$, where $\nu\colon G \to X$ is a 
$G$-compactification of $G$ such that $f(g)=F(\nu(g))$, and $F \in
\SUC(X)$. Clearly, $F$ is uniformly continuous because $X$ is
compact. Then $f \in \RUC(G)$ by Lemma \ref{l:right}. In order to
see that $f \in \LUC(G)$, choose $x_0:=\nu(e) \in X$ in the
definition of $\SUC$.
\end{proof}

\begin{remark} \label{r:abelian}
Recall that if $\mathcal{L}=\mathcal{R}$ then $G$ is said to be
SIN group. If $G$ is a SIN group then 
$X \in$ SUC for every compact $G$-space $X$. 
It follows that for a SIN group $G$ we have $\SUC(X)=\RUC(X)$ for
every, not necessarily compact, $G$-space $X$ and also
$\SUC(G)=\UC(G)=\RUC(G)$. For example this holds for
abelian, discrete, and compact groups.
\end{remark}

A special case of a SUC uniform $G$-space is obtained when the
uniform structure $\mu$ is defined by a $G$-invariant metric. If
$\mu$ is a metrizable uniformity then $(X,\mu)$ is uniformly
equicontinuous iff $\mu$ can be generated by a $G$-invariant
metric on $X$. A slightly sharper property is the local version: 
$\SUC_X \supset Eq_X$ (see Lemma \ref{l:LE}).

We say that a compactification of $G$ is {\em Roelcke\/} if the
corresponding algebra $\Acal$ is a $G$-subalgebra of $\UC(G)$, or
equivalently, if there exists a natural $G$-morphism $G^{\UC} \to
G^{\Acal}$.

\begin{lem} \label{l:SUC} We collect here the following properties of \emph{SUC}. 
\ben
\item $f(\SUC_X) \subset \SUC_Y$ for every uniformly continuous $G$-map $f \colon (X,\mu) \to (Y,\eta)$. 
\item
The class \emph{SUC} is closed under products, subsystems and quotients.
\item
 Let $\alpha\colon G \to Y$ be a Roelcke compactification.
Then $\alpha(G) \subset \SUC_Y$.
\item Let $X$ be a not necessarily compact $G$-space and $f \in
\SUC(X)$. Then for every $x_0 \in X$ and every $\eps >0$ there
exists a \nbd $U$ of $e$ such that
$$
|f(gux_0)-f(gx_0)| < \eps  \ \ \quad \forall \ (g,u) \in G \times
U.
$$
 \een
\end{lem}
\begin{proof} (1) and (2) are straightforward.

(3): \ Follows directly from (1) because the left action of $G$ on
itself is uniformity equicontinuous \wrt the left uniformity and
$\a\colon (G, \L) \to (Y, \mu_Y)$ is uniformly continuous for every
Roelcke compactification.

(4): \ There exist: a compact SUC $G$-space $Y$, a continuous
function $F\colon Y \to \R$ and a $G$-compactification $\nu\colon X \to Y$
such that $f=F \circ \nu$. Now our assertion follows from the fact
that $\nu(x_0) \in \SUC_Y$ for every $x_0 \in X$ 
(taking into account that $F$ is uniformly continuous).
\end{proof}

By Lemmas \ref{l:SUC}.2, \ \ref{l:pu} and the standard subdirect
product construction (see \cite[Proposition 2.9.2]{GM}) we can
derive the following facts.

\begin{lem} \label{p:qLE} \ 
\ben
\item
The collection $\SUC(X)$ is a $G$-subalgebra of $\RUC(X)$ for every
$G$-space $X$ and the corresponding Gelfand space
$|\SUC(X)|=X^{\SUC}$ with the canonical compactification $j \colon X \to
X^{\SUC}$ is the maximal SUC-compactification of~$X$.
\item
A compact $G$-system $X$ is \emph{SUC} if and only if $C(X)=\SUC(X)$.
\item For every $f \in \RUC(X)$ we have
$f \in \SUC(X)$ if and only if $X_f$ is \emph{SUC}.
\item
 $\SUC(G)$ is a {\it point-universal} closed $G$-subalgebra of
 $\RUC(G)$.
\item
The canonical compactification $j\colon G \to G^{\SUC}$ is always a
right topological semigroup compactification of $G$.
 \een
\end{lem}

We need also the following link between SUC functions and cyclic
$G$-spaces.

\begin{lem} \label{l:wle=luc}
Let $X$ be a $G$-space, $f \in \RUC(X)$ and $h \in X_f$. Then the
following are equivalent: \ben
\item
$h \in \LUC(G)$.
\item
$h \in \UC(G)$.
\item
$h$, as a point in the $G$-flow $Y:=X_f$, is in \emph{SUC}$_Y$. \een
\end{lem}
\begin{proof}
We know by Proposition \ref{f:X_f}.1 that $X_f \subset \RUC(G)$.
Thus (1) $\Longleftrightarrow$ (2).

For (1) $\Longleftrightarrow$ (3) observe that $$|h(gu) -h(g)| <
\eps \ \ \ \forall \ g \in G \quad \Longleftrightarrow \quad
|h(t_kgu) -h(t_kg)| < \eps \ \ \ \forall \ g \in G, \
k=1,\cdots,n$$ for arbitrary finite subset $\{t_1, \cdots, t_n\}$
of $G$.
\end{proof}

The following result shows that $\SUC(X)$ can be described by
internal terms for every compact $G$-space $X$.

\begin{prop} \label{p:SUC(X)}
Let $X$ be a compact $G$-space. The following are equivalent: \ben
\item  $f \in \SUC(X)$.
\item
For every $x_0 \in X$ and every $\eps >0$ there exists a \nbd $U$
of $e$ such that
$$
|f(gux_0)-f(gx_0)| < \eps  \ \ \quad \forall \ (g,u) \in G \times
U.
$$
 \een
\end{prop}
\begin{proof}
(1) $\Longrightarrow$ (2): \ Apply Lemma \ref{l:SUC}.4.

(1) $\Longleftarrow$ (2): \ By Lemma \ref{p:qLE}.3 we have to show
that the cyclic $G$-space $X_f$ is SUC. Fix an arbitrary element
$\om \in X_f$. According to Lemma \ref{l:wle=luc} it is equivalent
to verify that $\omega \in \LUC(G)$. The $G$-compactification map
$f_{\sharp}\colon X \to X_f$ is onto because $X$ is compact. Choose
$x_0 \in X$ such that $f_{\sharp}(x_0)=\om$.
Then $\omega(g)=f(gx_0)$ for every $g \in G$.
By assertion (2) for $x_0 \in X$ and $\eps$ we can pick a \nbd $U$
of $e$ such that
$$
|f(gux_0)-f(gx_0)| \leq \eps  \ \ \quad \forall \ (g,u) \in G
\times U.
$$
holds. Now we can finish the proof by observing that  
$$
|f(gux_0)-f(gx_0)| = |\om(gu) -\om(g)|  \leq \eps  \ \ \quad
\forall \ (g,u) \in G \times U.
$$
%
%
\end{proof}

The following result emphasizes the differences between $\RUC(G)$
and $\SUC(G)$.

\begin{thm} \label{t:SUC}
Let $\a=\a_{\mathcal{A}} \colon G \to S$ be a $G$-compactification of
$G$ such that
the corresponding left $G$-invariant subalgebra
$\mathcal{A}$ of $\RUC(G)$ is also right $G$-invariant. Consider
the following conditions: \ben
\item
$\mathcal{A} \subset \UC(G)$ (that is, $\a \colon G \to S$ is a Roelcke
compactification).
\item
The induced right action $S \times G \to S, \ (s,g) \mapsto s \a
(g)$ is jointly continuous.
\item
$\mathcal{A} \subset \SUC(G)$.
 \een
 Then
\bit
\item
[(a)] always, 1 $\Leftrightarrow$ 2 and 3 $\Rightarrow$ 1.
\item
[(b)] if, in addition, $S$ is a right topological semigroup and
$\a \colon G \to S$ is a right topological semigroup compactification of
$G$ then 1 $\Leftrightarrow$ 2 $\Leftrightarrow$ 3.
 \eit
\end{thm}
\begin{proof} (a)
1 $\Leftrightarrow$ 2:
By our assumption $\Acal$ is $G$-invariant with respect to
left and right translations (that is, the functions 
$(fg)(x):=f(gx)$ and $(gf)(x):=f(xg)$ lie in $\Acal$ for every
$f\in \Acal$ and $g \in G$). Then the corresponding (weak star
compact) Gelfand space $S:=X^{\Acal} \subset \Acal^*$, admits the
natural \emph{dual left and right actions} (see also Definition
\ref{d:dual} and Remark \ref{r:regular}) $\pi_l \colon G \times S \to S$
and $\pi_r \colon S \times G \to S$ such that $(g_1s)g_2=g_1(sg_2)$ for
every $(g_1,s,g_2) \in G \times S \times G$.
It is easy to see that this right action $S \times G \to S$
is jointly continuous if and only if $\mathcal{A} \subset \LUC(G)$.
On the other hand, since $\a \colon G \to S$ is a $G$-compactification of
left $G$-spaces we already have $\mathcal{A} \subset \RUC(G)$.

3 $\Rightarrow$ 1: By Lemma \ref{l:SUC-s-UC} we have $\SUC(G)
\subset \UC(G)$.

(b) We have to verify that 1 $\Rightarrow$ 3 provided that $\a \colon G
\to S$ is a right topological semigroup compactification of $G$.
The latter condition is equivalent to the fact that the system
$(G,S)$ is point universal (Lemma \ref{l:pu}) and thus for every
$x_0 \in S$ there exists a homomorphism of $G$-ambits $\phi:
(S,\a(e)) \to ({\cls}(Gx_0),x_0)$. By Lemma \ref{l:SUC} we
conclude that the point $x_0=\phi(\a(e)) \in S$ is a point of SUC
in the $G$-system ${\cls}(Gx_0)$ (and, hence, in $S$). Since $x_0$
is an arbitrary point in $S$ we get that $\SUC_S=S$ and hence $S$
is an SUC system.
Since every function $f \in \Acal$ on $G$ comes
from the compactification $\a \colon G \to S$ we conclude that
$\mathcal{A} \subset \SUC(G)$.
\end{proof}

\begin{cor} \label{c:SUC}
The $G$-compactification $j \colon G \to G^{\SUC}$ is a right topological 
semigroup compactification of $G$ such that the right action
$G^{\SUC} \times G \to G^{\SUC}$ is
also jointly continuous.
\end{cor}
\begin{proof}
Apply Proposition \ref{p:qLE}.5 and Theorem \ref{t:SUC}.
\end{proof}

\begin{cor} \label{p:CUSUC}
The following conditions are equivalent:
\ben
\item
$i \colon G \to G^{UC}$ is a right topological semigroup
compactification.
\item
$(G^{\UC},i(e))$ is a point universal $G$-system.
\item
$G^{\UC}$ is \emph{SUC}.
\item
$\SUC(G)=\UC(G)$.
\een
\end{cor}
\begin{proof}
 Apply assertion (b) of Theorem \ref{t:SUC} to
 $\mathcal{A}=\UC(G)$ 
taking into account Lemmas \ref{l:pu} and \ref{l:SUC-s-UC}. 
\end{proof}

Particularly interesting examples of groups $G$ with
$\SUC(G)=\UC(G)$ are the Polish groups $U(H)$ of all unitary
operators (Example \ref{e:unit}), and the group $S_\infty(\N)$
(Theorem \ref{env}). In both cases we actually have
$\SUC(G)=\UC(G)=\WAP(G)$. Note that these groups are not SIN (compare
Remark \ref{r:abelian}).  
%

For the next result see also Veech \cite[Section 5]{V}.

\begin{thm} \label{p:X_fin}
Let $f \in \RUC(X)$. The following conditions are equivalent: 
\ben
\item
$X_f \subset \UC(G)$.
\item
$(G,X_f)$ is \emph{SUC}.
\item
$f \in \SUC(X)$. 
\een
\end{thm}
\begin{proof}
1 $\Rightarrow$ 2: Let $h \in X_f$. By our assumption we have $h
\in \UC(G) \subset \LUC(G)$. Then by Lemma \ref{l:wle=luc}, $h$, as
a point in the $G$-flow $Y:=X_f$, is in SUC$_Y$. So $\SUC_Y=Y$.
This means that $(G,X_f)$ is SUC.

2 $\Rightarrow$ 3: Let $(G,X_f)$ be SUC. By Proposition
\ref{f:X_f}.3, the function $f \colon X \to \R$ comes from the
$G$-compactification $f_{\sharp} \colon X \to X_f$. By Definition
\ref{d:SUC}.3 this means that $f \in \SUC(X)$.

3 $\Rightarrow$ 1: Let $f \in \SUC(X)$. Then Lemma \ref{p:qLE}.3
says that $X_f$ is SUC. By Lemma \ref{l:wle=luc} we have $X_f
\subset \UC(G)$.
\end{proof}

\section{SUC, homogeneity and the epimorphism problem}

We say that a $G$-space $X$ is a \emph{coset $G$-space} if it is
$G$-isomorphic to the usual coset $G$-space $G/H$ where $H$ is a
closed subgroup of $G$ and $G/H$ is equipped with the quotient
topology. We say that a $G$-space $X$ is \emph{homogeneous} if for
every $x,y \in X$ there exists $g \in G$ such that $gx=y$. A
homogeneous $G$-space $X$ is a coset $G$-space if and only if the
orbit map $\tilde{x} \colon G \to X$ is open for some (equivalently,
every) $x \in X$. Furthermore, $\tilde{x} \colon G \to X$ is open iff it
is a quotient map. Recall that by a well known result of Effros
every homogeneous $G$-space with Polish $G$ and $X$ is necessarily
a coset $G$-space.

\begin{prop} \label{hom}
Let $X=G/H$ be a compact coset $G$-space. \ben
\item
If $X$ is \emph{SUC} then $X$ is equicontinuous (that is, almost
periodic).
\item
$\SUC(X)=\mathrm{AP}(X)$.
\een
\end{prop}
\begin{proof} (1): \ Indeed let $x_0H \in G/H$ and let
$\eps$ be an element of the uniform structure on the compact space
$X$. By Definition \ref{d:SUC} we can choose a \nbd $U$ of $e$
such that
$$
(gux_0H, gx_0H) \in \eps  \ \ \quad \forall \ (g,u) \in G \times
U.
$$
By the definition of coset space topology the set $O:=Ux_0H$ is a
\nbd of the point $x_0H$ in $G/H$. We obtain that $(gxH,gx_0H) \in
\eps$ whenever $xH \in O$. This proves that $x_0H$ is a point of
equicontinuity of $X=G/H$. Hence $X$ is AP.

(2): \ Every equicontinuous compact $G$-space is clearly SUC. This
implies that always, $\SUC(X) \supset \mathrm{AP}(X)$. Conversely, let $f
\in \SUC(X)$. This means that $f=\a \circ F$ for a
$G$-compactification $\a \colon X \to Y$ where $Y$ is SUC and $F \in
C(Y)=\SUC(Y)$. We can suppose that $\a$ is onto because $X$ is
compact. Then $\a$ is a quotient map. On the other hand $X$ is
a coset space $G/H$. It follows that the natural onto map $G \to Y$
is also a quotient map. Therefore, $Y$ is also a coset space of
$G$. Now we can apply (1). It follows that $Y$ is almost periodic.
Hence $f$ comes from an AP $G$-factor $Y$ of $X$. Thus, $f \in
\mathrm{AP}(X)$.
\end{proof}

Next we discuss a somewhat unexpected connection between SUC, free
topological $G$-groups and an \emph{epimorphism problem}.
Uspenskij has shown in \cite{Us-epic} that in the category of
Hausdorff topological groups epimorphisms need not have a dense
range. This answers a longstanding problem by K. Hofmann. Pestov
established \cite{Pe-epic,pest-wh} that the question completely
depends on the \emph{free topological $G$-groups} $F_G(X)$ of a
$G$-space $X$ in the sense of Megrelishvili
\cite{Me-free-G-group}. More precisely, the inclusion $i \colon H
\hookrightarrow G$ of topological groups is epimorphism iff the
free topological $G$-group $F_G(X)$ of the coset $G$-space
$X:=G/H$ is trivial. Triviality means, `as trivial as possible',
isomorphic to the cyclic discrete group.

For a $G$-space $X$ and points $x_1, x_2 \in X$ we write $x_1
\stackrel{Aut} {\sim} x_2$ if these two points have the same image
under the canonical $G$-map $X \to F_G(X)$.
If $d$ is a 
bounded compatible $G$-invariant metric
on a $G$-space $X$ then 
$(X,d)$ is isometrically $G$-linearizable using Arens-Eells embedding. 
 Therefore, in this case $x_1
\stackrel{Aut} {\sim} x_2$ iff $x_1=x_2$.

\begin{thm} \label{t:epic}
Let $H$ be a closed subgroup of $G$.
\ben
\item
If $x_1 \stackrel{Aut} {\sim} x_2$ for $x_1, x_2$ in the $G$-space
$X:=G/H$ then $x_1 \stackrel{\SUC} {\sim} x_2$.
\item
If the inclusion $H \hookrightarrow G$ is an epimorphism then the
coset $G$-space $G/H$ is \emph{SUC}-trivial.   
\een
\end{thm}
\begin{proof} (1)
Assuming the contrary let $f\colon G/H \to \R$ be SUC function which
separates the points $x_1:=a_1H$ and $x_2:=a_2H$. Then the bounded 
$G$-invariant pseudometric $\rho_f$ on $G/H$ defined by
$\rho_f(xH,yH):=\sup_{g \in G} |f(gxH)-f(gyH)|$ also separates
these points. We show that $\rho_f$ is continuous. Indeed let
$\eps >0$ and $x_0H \in G/H$. By virtue of Lemma \ref{l:SUC}.4 we
can choose a \nbd $U$ of $e$ such that
$$
|f(gux_0H)-f(gx_0H)| < \eps  \ \ \quad \forall \ (g,u) \in G
\times U.
$$
By the definition of coset space topology the set $O:=Ux_0H$ is a
\nbd of the point $x_0H$ in $G/H$. We obtain that $\rho_f(xH,
x_0H) < \eps$ whenever $xH \in O$. This proves the continuity of
$\rho_f$.

Consider the associated metric space $(Y,d)$ and the canonical
distance preserving onto $G$-map $X \to Y, x \mapsto [x]$. The
metric $d$ on $Y$ (defined by $d([x],[y]):=\rho_f(x,y)$) is
$G$-invariant. 
 Since $d([x_1], [x_2]) >0$ we conclude 
that $x_1$ and $x_2$ have different images in $F_G(X)$ (see the discussion above). This
contradicts the assumption $x_1 \stackrel{Aut} {\sim} x_2$.

(2) Assume that $G/H$ is not SUC-trivial. By Assertion (1) we get
that the free topological $G$-group $F_G(G/H)$ of $G/H$ is not
trivial. Therefore by the above mentioned result of Pestov
\cite{Pe-epic} we can conclude that the inclusion $H
\hookrightarrow G$ is not an epimorphism.
\end{proof}

\begin{remark} \label{r:epic} \ 
\ben
\item
The converse to Theorem \ref{t:epic}.2 is not true (take
$G:=H_+[0,1]$, \ $H:=\{e\}$ and apply Theorem \ref{interval}).
\item
As a corollary of Theorem \ref{t:epic}.2 one can get several
examples of SUC-trivial (compact) $G$-spaces. For example, by
\cite{Me-free-G-group} the free topological $G$-group $F_G(X)$ of
$X:=G/H$ with $G:=H(\mathbb{T}), \ H:= St(z)$
(where $z \in \mathbb{T}$ is an
arbitrary point of the circle $\mathbb{T}$) is trivial.
In fact it is easy to see that the same is true for the smaller
group $G:=H_+(\mathbb{T})$
(and the subgroup $H:= St(z)$) (cf., Proposition \ref{hom}).
\item
It is a well known result by Nummela \cite{Num} that if $G$ is a
SIN group then the inclusion of a closed proper subgroup $H
\hookrightarrow G$ is not an epimorphism. This result easily
follows from Theorem \ref{t:epic}.2. Indeed if $G$ is SIN then by
Remark \ref{r:abelian} for the coset $G$-space $G/H$ we have
$\SUC(G/H)=\RUC(G/H)$. Hence if $G/H$ is SUC-trivial then
necessarily $H=G$ because $\RUC(G/H)$ is non-trivial for every
closed proper subgroup $H$ of $G$.
\een
\end{remark}

\section{Representations of groups and $G$-spaces on Banach spaces}
\label{s:repr}

For a real normed space $V$ denote by $B_V$ its closed unit ball
$\{v \in V:  ||v|| \leq 1 \}$. Denote by $\Iso(V)$ the topological
group of all linear surjective isometries $V \to V$ endowed with
the \emph{strong operator topology}. This is just the topology of
pointwise convergence inherited from $V^V$. Let $V^*$ be the dual
Banach space of $V$ and
$$
\langle\ ,\ \rangle \,  \colon V \times V^* \to \R,
\quad \ (v,\psi) \mapsto \, \langle v,\psi \rangle = \psi(v)
$$
is the canonical (always continuous) bilinear mapping.

A \emph{representation} (\emph{co-representation}) of a
topological group $G$ on a normed space $V$ is a homomorphism
(resp. co-homomorphism) $h \colon G \to \Iso(V)$.
Sometimes it is more convenient to describe a representation
(co-representation) by the corresponding left (resp. right) linear
isometric actions $\pi_h \colon G \times V \to V, \ (g,v) \mapsto
gv=h(v)(g)$ (resp., $V \times G \to V, \ (v,g) \mapsto
vg=h(v)(g)$).
The (co)representation $h$ is continuous if and only if the action
$\pi_h$ is continuous.

\begin{remark} \label{r:co}
Many results formulated for co-representations remain true also
for representations (and vice versa) taking into account the
following simple fact: for every representation
(co-representation) $h$ there exists an associated
co-representation (representation) $h^{op} \colon G \to \Iso(V), \ g
\mapsto h(g^{-1}).$
\end{remark}

\begin{defn} \label{d:dual}
Let $\pi \colon G \times V \to V$ be a continuous left action of
$G$ on $V$ by linear operators. The \emph{adjoint (or, dual)
right} action $\pi^* \colon  V^* \times G \to V^*$ is defined by $\psi
g(v):=\psi(gv)$. The corresponding \emph{adjoint (dual) left}
action is $\pi^* \colon G \times V^* \to V^*$, where
$g\psi(v):=\psi(g^{-1}v)$. Similarly, if $\pi \colon V \times G \to V$
is a continuous linear \emph{right} action of $G$ on $V$ (e.g.,
induced by some co-representation), then the corresponding
\emph{adjoint (dual) action} $\pi^* \colon G \times V^* \to V^*$ is
defined by $g\psi(v):=\psi(vg)$.
\end{defn}

The main question considered in \cite{Mefr} was whether the dual
action $\pi^*$ of $G$ on $V^*$
is jointly continuous with respect to the norm topology on $V^*$.
When this is the case we say that the action $\pi$ (and, also the
corresponding representation $h \colon G \to \Iso(V)$, when $\pi$ is an
action by linear isometries) is {\it adjoint continuous}.
This name was suggested by V. Uspenskij.

\begin{remark}
In general, not every continuous representation is adjoint
continuous (see for example \cite{Mefr}). A standard example is
the representation of the circle group
$G:=\mathbb{T}$ on $V:=C(\mathbb{T})$ by
translations. Here the Banach space $V$ is separable but with
``bad" geometry. The absence of adjoint continuity may happen even
for relatively ``good" (for instance, for separable
Radon-Nikod\'ym) Banach spaces like $V:=l_1$. Indeed if we
consider the symmetric group $G:=S_{\infty}$, naturally embedded into
$\Iso(V)$ (endowed with the strong operator topology) as the group
of ``permutation of coordinates" operators, then
the dual action of $G$ on $l_1^*=l_{\infty}$ is not continuous
(see \cite{Meop}).
\end{remark}

It turns out that the situation in that respect is the best
possible for the important class ${\mathcal Asp}$ of Asplund
Banach spaces. The investigation of this class and the closely
related \emph{Radon-Nikod\'ym property} is among
the main themes in
Banach space theory. Recall that a Banach space $V$ is an {\em
Asplund\/} space if the dual of every separable linear subspace is
separable, iff every bounded subset $A$ of the dual $V^*$ is
(weak${}^*$,norm)-\emph{fragmented}, iff $V^*$ has the
Radon-Nikod\'ym property. Reflexive spaces and spaces of the type
$c_0(\Gamma)$ are Asplund. For more details cf. \cite{Bo,Fa}. For
the readers convenience we recall also the definition of
fragmentability.

\begin{defn} \label{d:fragm}
(Jane and Rogers \cite{JR}) Let $(X,\tau)$ be a topological space
and $\rho$ be a metric on the set $X$. Then $X$ is said to be
$(\tau,\rho)$-\emph{fragmented} if for every nonempty $A \subset
X$ and every $\eps >0$ there exists a $\tau$-open subspace $O$ of
$X$ such that $O \cap A$ is nonempty and $\eps$-small in
$(X,\rho)$.
\end{defn}

Namioka's Joint continuity theorem implies that every weakly
compact
set in a Banach space is norm fragmented. This explains why
every reflexive space is Asplund.


\begin{thm} \label{t:adj}
\cite[Corollary 6.9]{Mefr} Let $V$ be an Asplund Banach space. If
a (not necessarily isometric) linear action $\pi \colon G \times V \to
V$
is continuous then the dual right action $\pi^* \colon V^* \times G \to
V^*$ is also continuous.
\end{thm}

Certainly, this result remains true for dual left actions $\pi^*:
G \times V^* \to V^*$, where $g\psi(v):=\psi(g^{-1}v)$, as well as
for dual actions defined by a right action $\pi \colon V \times G \to V$.
The obvious reason is the continuity of the map $G \to G, \ g
\mapsto g^{-1}$.

\sk

The following definition provides a flow version of the group
representation definitions discussed above. It differs from the
usual notion of $G$-\emph{linearization} in that here we represent
the phase space of the flow as a subset of the dual space $V^*$
(with respect to the dual action and \emph{weak star topology})
rather than as a subset of $V$.

\begin{defn} \label{d:repr}
\cite{Menz} Let $X$ be a $G$-space. A continuous (\emph{proper})
\emph{representation} of $(G,X)$ on a Banach space $V$ is a pair
$$(h,\a) \colon G\times X \rightrightarrows {\Iso}(V)\times B^*$$ where
$h \colon G \to \Iso(V)$ is a strongly continuous co-homomorphism and
$\a \colon X \to B^*$ is a weak star continuous $G$-map (resp.
\emph{embedding}) with respect to the {\em dual action\/} $G
\times V^* \to V^*, \ (g\varphi)(v):=\varphi(h(g)(v))$. Here $B^*$
is the weak star compact unit ball of the dual space $V^*$.

Alternatively, one can define a representation in such a way that
$h$ is a \emph{homomorphism} and the dual action $G \times V^* \to
V^*$ is defined by $(g\varphi)(v):=\varphi(h(g^{-1})(v))$.
\end{defn}

\begin{remark} \label{r:regular}
Let $X$ be a $G$-space and let $\Acal$ be a Banach (closed,
unital) subalgebra of $C(X)$. Associated with $\Acal$ we have the
canonical $\Acal$-compactification $\nu_{\Acal} \colon X \to X^{\Acal}$
of $X$, where $X^{\Acal}=|\Acal|$ is the Gelfand space of $\Acal$.
Here $X^{\Acal}$ is canonically embedded into the weak star
compact unit ball $B^*$ of the dual space $\Acal^*$. If $\Acal$ is
{\it $G$-invariant} (that is, the function $(fg)(x):=f(gx)$ lies
in $\Acal$ for every $f\in \Acal$ and $g \in G$) then $X^{\Acal}$
admits the natural adjoint action $G \times X^{\Acal} \to
X^{\Acal}$ with the property that all translations $\breve{g}:
X^{\Acal} \to X^{\Acal}$ are continuous and such that $\a_{\Acal}:
X \to X^{\Acal} \subset B^*$ is $G$-equivariant. We obtain in this
way a representation (where $h$ is not necessarily continuous)
$$
(h,\a_{\Acal}) \colon (G,X) \rightrightarrows (\Iso(\Acal), B^*)
$$
on the Banach space $\Acal$, where $h(g)(f):=fg$ (and
$\a_{\Acal}(x)(f):=f(x)$). We call it the {\it canonical} (or,
\emph{regular}) {\it $\Acal$-representation} of $(G,X)$. It is
continuous iff $\Acal \subset \RUC(X)$ (see for example \cite[Fact
2.2]{Menz} and \cite[Fact 7.2]{opit}). The regular
$\RUC(X)$-representation leads to the maximal $G$-compactification
$X \to \beta_GX$ of $X$. It is proper if and only if $X$ is
$G$-compactifiable.
\end{remark}

The following observation due to Teleman is well known (see also
\cite{pest-wh}).

\begin{fact}  \label{Teleman}
\emph{(Teleman \cite{Te})} Every topological group can be embedded
into $\Iso(V)$ for some Banach space $V$.
\end{fact}
\begin{proof}
It is well known that $\RUC(G)$ determines the topology of $G$.
Hence the regular $V:=\RUC(G)$-representation $(h,\a) \colon (G,G)
\rightrightarrows (\Iso(V), B^*)$ is proper. That is, the map $\a$
is an embedding. In fact it is easy to see that the
co-homomorphism $h$ is an embedding of topological spaces. The
representation $h^{op} \colon G \to \Iso(V), g \mapsto h(g^{-1})$ is then
a topological group embedding.
\end{proof}

\begin{defn} \label{d:GrRep}
Let $\mathcal{K}$ be a ``well behaved" subclass of the class
${\mathcal Ban}$ of all Banach spaces. Typical and important
particular cases for such ${\mathcal K}$ are: ${\mathcal Hilb}, \
{\mathcal Ref}$ or ${\mathcal Asp}$, the classes of Hilbert,
reflexive or Asplund Banach spaces respectively.
\ben
\item
A topological group $G$ is $\mathcal{K}$-\emph{representable} if
there exists a (co)representation $h \colon G \to \Iso(V)$ for some $V
\in \mathcal{K}$ such that $h$ is topologically faithful (that is,
an embedding). Notation: $G \in \mathcal{K}_r$.
\item
In the opposite direction, we say that $G$ is
$\mathcal{K}$-\emph{trivial} if every continuous
$\mathcal{K}$-representation
(or, equivalently, \emph{co-representation}) $h \colon G \to \Iso(V)$ is
trivial.
%
\item
We say that a topological group $G$ is SUC-\emph{trivial} if
$\SUC(G)=\{constants\}$. Analogously can be defined
WAP-\emph{trivial} groups. $G$ is WAP-trivial iff $G$ is
\emph{reflexively trivial} (${\mathcal Ref}$-trivial in the sense
of Definition \ref{d:GrRep}.2). Similarly, $Asp(G)=\{constants\}$
iff $G$ is ${\mathcal Asp}$-trivial. These equivalences follow for
instance from Theorem \ref{t:matrix} below.
\een
\end{defn}

\begin{remark} \label{r:repr} \ 
\ben
\item
By Teleman's theorem (Fact \ref{Teleman}) every topological group
is ``Banach representable". Hence,
$\{Topological \ Groups\}={\mathcal Ban}_r$.
\item
$\{Topological \ Groups\}={\mathcal Ban}_r \supset {\mathcal
Asp}_r \supset {\mathcal Ref}_r \supset {\mathcal Hilb}_r.$
\item \ By
Herer and Christensen \cite{HC} (see also Banasczyk \cite{Ba})
abelian (even monothetic) groups can be ${\mathcal Hilb}$-trivial.
Note also that $c_0  \notin {\mathcal Hilb}_r$, \cite{Meist}.

\item \ \cite{Meist} 
The additive group $L_4[0,1]$ is reflexively but not Hilbert
representable.
\item \ \cite{Merup} \ $H_+[0,1] \notin {\mathcal Ref}_r$.
It was shown in \cite{Merup} that every weakly almost periodic
function on the topological group $G:=H_+[0,1]$ is constant and
that $G$ is ${\mathcal Ref}$-trivial. By Pestov's observation (see
\cite[Corollary 1.4]{Pe-new} and Lemma \ref{l:observation}) the
same is true for the group $\Iso(\U_1)$.
\item \ Theorem \ref{t:conclusions}.3 shows that $H_+[0,1]$ is even
${\mathcal Asp}$-trivial. In fact we show that every ``adjoint
continuous" representation of that group is trivial (Theorem
\ref{t:conclusions}.2). This result was obtained also by Uspenskij
(unpublished). Furthermore, we prove a stronger result by showing
that $H_+[0,1]$ (and also $\Iso(\U_1)$) are SUC-trivial. \een
\end{remark}

\begin{problem} \label{q:asp} (See also \cite{opit} and \cite{Merup})
\ben
\item
Distinguish ${\mathcal Asp}_r$ and ${\mathcal Ref}_r$ by
 finding $G \in {\mathcal Asp}_r$ such that $G \notin {\mathcal
 Ref}_r$.
\item
Find an abelian $G \notin {\mathcal Ref}_r$.
 \een
\end{problem}

\sk

Now we turn to the ``well behaved actions". Recall the dynamical
versions of Eberlein and Radon-Nikod\'ym compact spaces.

\begin{defn} \label{d-RN}
\cite{GM, Menz} Let $X$ be a $G$-space.
\begin{enumerate}
\item
$(G,X)$ is a {\em Radon-Nikod\'ym system\/} (RN for short) if
there exists a proper representation of $(G,X)$ on an Asplund
Banach space $V$. If we can choose $V$ to be reflexive, then
$(G,X)$ is called an {\em Eberlein system}. The classes of
Radon-Nikod\'ym and Eberlein compact systems will be
denoted by RN and Eb respectively.
\item
$(G,X)$ is called an RN-{\em approximable\/} system
($\rm{RN_{app}}$) if it can be represented as a subdirect product
of RN systems.
\end{enumerate}
\end{defn}

Note that compact spaces which are not Eberlein are necessarily
non-metrizable, while even for $G:=\Z$, there are many natural
{\it metric} compact $G$-systems which are not RN.

\begin{defn} \label{d:adj} \ 
\ben
\item
A representation $(h,\a)$ of a $G$-space $X$ on $V$ is
\emph{adjoint continuous} if the dual action $G \times V^* \to
V^*$ is also continuous (or, equivalently, if the group
corepresentation $h \colon G \to \Iso(V)$ is adjoint continuous).
\item
Denote by ${\mathcal Adj}$ the class of compact $G$-systems which
admit a proper adjoint continuous representation on some Banach
space $V$. Theorem \ref{t:adj} implies that RN $\subset {\mathcal
Adj}$.
\item
 Denote by ${\adj}(G)$ the collection of functions on $G$
which come from a compact $G$-space $X$ such that $(G,X)$ is in
the class ${\mathcal Adj}$. In fact this means that $f$ can be
represented as a generalized matrix coefficient (see Section
\ref{s:mat-SUC}) of some adjoint continuous representation of $G$.
\een
\end{defn}

\begin{prop} \label{mat-qLE}
$\Asp(G) \subset {\rm{adj}}(G)$ for every topological group $G$.
\end{prop}
\begin{proof} By \cite[Theorem 7.11]{Menz}
(or Proposition \ref{p:asp-cond}) $f \in \Asp(G)$ iff $f$ comes
from a $G$-compactification $G \to X$ of $G$ with $X \in$ RN. Now
observe (as in Definition \ref{d:adj}.2) that RN $\subset
{\mathcal Adj}$ by Theorem \ref{t:adj}.
\end{proof}

\section{Dynamical complexity of functions}
\label{s:hierarchy}

In this section we introduce a hierarchy of dynamical complexity
of functions on a topological group $G$ which reflects the
complexity of the $G$-systems from which they come. Our main tool
is the cyclic $G$-system $X_f$ corresponding to a function $f \colon X
\to \R$. Recall that when $X:=G$, the space $X_f$ is the pointwise
closure of the orbit $Gf$ in $\RUC(G)$. The topological nature of
$X_f$ in the Banach space $\RUC(G)$ relates to the dynamical
complexity of $f$ and leads to a natural hierarchy of complexity
(see Theorem \ref{r:hierarchy} below). In particular we will
examine the role that SUC functions play in this hierarchy.

\bsk \emph{Periodic orbits and the profinite compactification}
\sk

The most elementary dynamical system is a finite (periodic) orbit.
It corresponds to a clopen subgroup $H < G$ of finite index. These
subgroups form a directed set and the corresponding compact
inverse limit $G$-system
$$
X^{PF}=\underset{\leftarrow}{\lim}\ G/H
$$
is the
{\it profinite} compactification of $G$.

\bsk
\emph{Almost Periodic functions and the Bohr
compactification} \sk

The weaker requirement that $X_f$ be norm compact in $\RUC(G)$
leads to the well known definition of almost periodicity. A
function $f \in C(X)$ on a $G$-space $X$ is {\em almost
periodic\/} if the orbit $fG:=\{fg\}_{g \in G}$ forms a precompact
subset of the Banach space $C(X)$. The collection $\mathrm{AP}(X)$ of AP
functions is a $G$-subalgebra in $\RUC(X)$. The universal almost
periodic compactification of $X$ is the Gelfand space $X^{\mathrm{AP}}$ of
the algebra $\mathrm{AP}(X)$. When $X$ is compact this is the classical
{\em maximal equicontinuous} factor of the system $X$.  A compact
$G$-space $X$ is equicontinuous iff $X$ is almost periodic (AP),
that is, iff $C(X)=\mathrm{AP}(X)$. For a $G$-space $X$ the collection
$\mathrm{AP}(X)$ is the set of all functions which come from equicontinuous
(AP) $G$-compactifications.

For every topological group $G$, treated as a $G$-space, the
corresponding universal AP compactification is the well known
\emph{Bohr compactification} $b \colon G \to bG$, where $bG$ is a
compact topological group.

\begin{thm} \label{l:X_fandAP}
Let $X$ be a $G$-space. For $f \in \RUC(X)$ the
following conditions are equivalent:
\ben
\item
$f \in \mathrm{AP}(X)$.
\item
$(G,X_f)$ is equicontinuous.
\item
$X_f$ is norm compact in $\RUC(G)$.
\een
\end{thm}
\begin{proof}
(1) $\Longleftrightarrow$ (2): $f \in \mathrm{AP}(X)$ iff the cyclic
algebra $\Acal_f$ (which, by Proposition \ref{f:X_f}, generates
the compactification $X \to X_f$) is a subalgebra of $\mathrm{AP}(X)$.

(2) $\Longleftrightarrow$ (3): It is easy to see that the
$G$-space $X_f$ is equicontinuous iff the norm and pointwise
topologies coincide on $X_f \subset \RUC(G).$
\end{proof}

\bsk
\emph{Weakly Almost Periodic functions}
\sk

A function $f \in C(X)$ on a $G$-space $X$ is called {\em weakly
almost periodic\/} (WAP for short; notation: $f \in \WAP(X)$) if
the orbit $fG:=\{fg\}_{g \in G}$ forms a weakly precompact subset
of $C(X)$. A compact $G$-space $X$ is said to be {\em weakly
almost periodic\/} \cite{E0} if $C(X)=\WAP(X)$.
 For a $G$-space $X$
the collection $\WAP(X)$ is the set of all functions which come
from WAP $G$-compactifications.
%
The universal WAP $G$-compactification $X \to X^{\WAP}$ is well
defined. The algebra $\WAP(G)$ is a point-universal $G$-algebra
containing $\mathrm{AP}(G)$. The compactification $G \to G^{\WAP}$ (for $X:=G$) is the universal semitopological semigroup
compactification of $G$.

A compact $G$-space $X$ is WAP iff it admits sufficiently many
representations on reflexive Banach spaces \cite{Menz}.
Furthermore if $X$ is a metric compact $G$-space then $X$ is WAP
iff $X$ admits a proper representation on a reflexive Banach
space. That is, iff $X$ is an Eberlein $G$-space.

\begin{thm} \label{l:WAPandX_f}
Let $X$ be a $G$-space. For $f \in \RUC(X)$ the following
conditions are equivalent: \ben
\item
$f \in \WAP(X).$
\item
$(G,X_f)$ is \emph{WAP}.
\item
$X_f$ is weak compact in $\RUC(G)$.
\item
$(G,X_f)$ is Eberlein (i.e., reflexively representable).
\een
\end{thm}

\begin{proof}
(1) $\Longleftrightarrow$ (2): $f \in \WAP(X)$ iff the algebra
$\Acal_f$ is a subalgebra of $\WAP(X)$.

(2) $\Longrightarrow$ (3): Let $F_e \colon X_f \to \R, \
F_e(\omega)=\om(e)$ be as in the definition of $X_f$. Consider the
weak closure $Y:={\cls}_w(F_eG)$ of the orbit $F_eG$.
Then $Y$ is weakly compact in $C(X_f)$ because $F_e \in
C(X_f)=\WAP(X_f)$ is weakly almost periodic. If $\om_1$ and $\om_2$
are distinct elements of $X_f$ then $(F_eg)(\om_1)=\om_1(g) \neq
\om_2(g)=(F_eg)(\om_2)$ for some $g \in G$. This means that the
separately continuous evaluation map $Y \times X_f \to \R$
separates points of $X_f$. Now $X_f$ can be treated as a pointwise
compact bounded subset in $C(Y)$. Hence by Grothendieck's well
known theorem \cite{Groth} we get that $X_f$ is weakly compact in
$C(Y)$. Since $G \to Y, \ g \mapsto gF_e$ is a
$G$-compactification of $G$, we have a natural embedding of Banach
algebras $j \colon C(Y) \hookrightarrow \RUC(G)$. It follows that
$X_f=j(X_f)$ is also weakly compact as a subset of $\RUC(G)$.

(3) $\Longrightarrow$ (4): The isometric action $G \times \RUC(G)
\to \RUC(G), \ (g,f) \mapsto gf$ induces a representation $h \colon G \to
\Iso(\RUC(G))$. If the $G$-subset $X_f$ is weakly compact in
$\RUC(G)$ then one can apply Theorem 4.11 (namely, the equivalence
between (i) and (ii)) of \cite{Menz} which guarantees that the
$G$-space $X_f$ is Eberlein.

(4) $\Longrightarrow$ (1): $f \in \WAP(X)$ because it comes from
$(G,X_f)$ (Proposition \ref{f:X_f}) which is WAP 
(being reflexively representable). 
\end{proof}

\bsk \emph{Asplund functions, ``sensitivity to initial conditions"
and Banach representations}
\sk

The following definition of ``sensitivity to initial conditions"
is essential in several definitions of chaos in dynamical systems,
mostly for $G:=\Z$ or $\R$ actions on metric spaces (see for
instance papers of Guckenheimer, Auslander and Yorke, Devaney,
Glasner and Weiss).

\begin{defn} \label{d:sens} \cite{GM}
Let $(X,\mu)$ be a uniform $G$-space.
\ben
\item
We say that $X$ is \emph{sensitive to initial conditions}
(or just \emph{sensitive}) if there
exists an $\eps \in \mu$ such that for
every nonempty open subset $O$
of $X$ the set $gO$ is not $\eps$-small for some $g \in G$.
Otherwise, $X$ is \emph{non-sensitive} (for short: NS).
\item
$X$ is \emph{Hereditarily Non Sensitive} (HNS) if every closed
$G$-subspace of $X$ is NS. \een
\end{defn}

Denote by HNS the class of all compact HNS systems. The following
result says that a compact $G$-system $X$ is HNS iff $(G,X)$
admits sufficiently many representations on Asplund spaces.

\begin{thm} \label{f:HNS} \cite{GM}
%
\ben
\item
\ HNS = RN$_{app}$.
\item
If $X$ is a compact metric $G$-space then $X$ is HNS iff $X$ is RN
(that is, Asplund representable). \een
\end{thm}

A function $f \colon X \to \R$ on a $G$-space $X$ is \emph{Asplund}
(notation: $f \in \Asp(X)$) \cite{Menz} if it satisfies one of the
following equivalent conditions.

\begin{prop} \label{p:asp-cond}
Let $f \colon X \to \R$ be a function on a $G$-space $X$. The following
conditions are equivalent: \ben
\item
$f$ comes from a $G$-compactification $\nu \colon X \to Y$ where $(G,Y)$
is HNS.
\item
$f$ comes from a $G$-compactification $\nu \colon X \to Y$ where $(G,Y)$
is RN.
\item
$f$ comes from a $G$-compactification $\nu \colon X \to Y$ and a
function $F \colon Y \to \R$ where the pseudometric space
$(Y,\rho_{H,F})$ with
$$
\rho_{H,F}(x,x')=\sup_{h\in H}|F(h x) - F(h x')|.
$$
is separable for every countable (equivalently, second countable)
subgroup $H \subset G$. \een
\end{prop}

The collection $\Asp(X)$ is always a $G$-subalgebra of $\RUC(X)$. It
defines the maximal HNS-compactification $X \to X^{\Asp}=|\Asp(X)|$
of $X$. For every topological group $G$ the algebra $\Asp(G)$ (as
usual, $X:=G$ is a left $G$-space) is point-universal.

\begin{thm} \label{l:HNSandX_f}
Let $X$ be a $G$-space. For every $f \in \RUC(X)$
the following conditions are equivalent:
\ben
\item
$f \in \Asp(X)$.
\item
$(G,X_f)$ is RN.
\item
$X_f$ is norm fragmented in $\RUC(G)$.
\een
\end{thm}

\begin{proof} If $X$ is compact then the proof
follows directly from \cite[Theorem 9.12]{GM}. Now observe that
one can reduce the case of an arbitrary $G$-space $X$ to the case
of a compact $G$-space $X_f$ by considering the cyclic $G$-system
$(X_f)_{F_e}$ (defined for $X:=X_f$ and $f:=F_e$) which can be
naturally identified with $X_f$.
\end{proof}

Explicitly the fragmentability of $X_f$ means that for every $\eps
>0$ every nonempty (closed) subset $A$ of $X_f \subset \R^G$
contains a relatively open (in the pointwise topology) nonempty
subset $O \cap A$ which is $\eps$-small in the Banach space
$\RUC(G)$.

As we already mentioned every weakly compact set is norm
fragmented so that $\WAP(X) \subset \Asp(X)$ for every $G$-space
$X$. In particular, $\WAP(G) \subset \Asp(G)$.

\bsk \emph{Locally equicontinuous functions} \sk

During the last decade various conditions weakening the classical
notion of equicontinuity were introduced and studied (see e.g.
\cite{GW2}, \cite{AAB1}, \cite{AAB2}, \cite{AG}). The following
definition first appears in a paper of Glasner and Weiss
\cite{GW2}.

\begin{defn} \label{d:1} \cite{GW2}
Let $(X,\mu)$ be a uniform $G$-space.
%
A point $x_0 \in X$ is a point of {\it local equicontinuity}
(notation: $x_0 \in \LE_X$) if $x_0$ is a point of equicontinuity in the uniform
$G$-subspace ${\cls}(Gx_0)$. We have $x_0 \in \LE_X$ iff $x_0 \in
\LE_Y$ iff $x_0 \in Eq_Y$ where $Y$ is the orbit $Gx_0$ of $x_0$
(see Lemma \ref{l:LE}.1). If $\LE_X=X$, then $X$ is {\it locally
equicontinuous} (LE).
\end{defn}

\begin{lem} \label{l:LE} \  
\ben
\item
Let $Y$ be a dense $G$-subspace of $(X,\mu)$ and $y_0 \in Y$. Then
$y_0 \in Eq_X$ if and only if $y_0 \in Eq_Y$.
\item
$\SUC_X \supset \LE_X \supset Eq_X$.
\item
$\SUC(X) \supset \LE(X)$.
\een
\end{lem}
\begin{proof}
(1)
 Let $\eps \in \mu$. There exists $\delta \in \mu$ such that
$\delta$ is a closed subset of $X\times X$ and $\delta \subset
\eps$. If $y_0 \in Eq_Y$ there exists an open set $U$ in $X$ such
that $y_0 \in U$ and $(gy,gy_0) \in \delta$ for all $y \in U \cap
Y$ and $g \in G$. Since $Y$ is dense in $X$ and $U$ is open we
have $U \subset {\cls}(U\cap Y)$. Since every $g$-translation $X
\to X, \ x \mapsto gx$ is continuous and $\delta$ is closed we get
$(gx,gy_0) \in \delta \subset \eps$ for every $g \in G$ and $x \in
U$.

(2): \ Let $x_0 \in \LE_X$. For every $\eps \in \mu$ there exists a
neighborhood $O(x_0)$ such that $(gx,gx_0)$ is $\eps$-small for
every $x \in O(x_0)$ and $g \in G$. Choose a neighborhood $U(e)$
such that $Ux_0 \subset O$. Then $(gux_0,gx_0)$ is $\eps$-small,
too. This proves the non-trivial part $\SUC_X \supset \LE_X$.

(3): Directly follows from (2).
\end{proof}

The collection $\LE(X)$ forms a $G$-subalgebra of $\RUC(X)$. Always,
 $Asp(X) \subset \LE(X)$. The algebra $\LE(X)$ defines
 the maximal LE-compactification $X \to X^{\LE}$ of $X$.
For every topological group $G$ the algebra $\LE(G)$ is
point-universal \cite{GM}.

\begin{thm} \label{l:LEand} \cite{GM} \
Let $X$ be a $G$-space. For $f \in \RUC(X)$ the following
conditions are equivalent: \ben
\item
$f \in \LE(X)$.
\item
$(G,X_f)$ is \emph{LE}.
\item
$X_f$ is \emph{orbitwise light} in $\RUC(G)$ (that is, for every
function $\psi \in X_f$ the pointwise and norm topologies coincide
on the orbit $G\psi$).
\een
\end{thm}
\begin{proof} (1) $\Longleftrightarrow$ (2) Directly follows
from \cite[Theorem 5.15.1]{GM}. On the other hand by \cite[Lemma
5.18]{GM} we have (2) $\Longleftrightarrow$ (3).
\end{proof}


\bsk \emph{The dynamical hierarchy} \sk


\begin{thm} \label{l:inclusions} \ 
\begin{enumerate}
\item
Let $X$ be a (not necessarily compact) $G$-space. We have the
following inclusions of $G$-subalgebras:
$$
\RUC(X) \supset \SUC(X) \supset \LE(X) \supset \Asp(X) \supset \WAP(X) \supset \mathrm{AP}(X)
$$
and the corresponding chain of $G$-factor maps
$$
\beta_G X \to X^{\SUC} \to X^{ LE} \to X^{Asp} \to X^{WAP} \to X^{
AP}.
$$
\item
For every topological group $G$ we have the following inclusions of $G$-subalgebras:
$$
\RUC(G) \supset \UC(G) \supset \SUC(G)
\supset \LE(G) \supset \Asp(G) \supset \WAP(G) \supset \mathrm{AP}(G)
$$ and
the corresponding chain of $G$-factor maps
$$
G^{\RUC} \to G^{\UC} \to G^{\SUC}
\to G^{\LE} \to G^{\Asp} \to G^{\WAP} \to G^{\mathrm{AP}}. 
$$
\end{enumerate}
\end{thm}

\begin{proof}
For the assertions concerning $\SUC(X)$ and $\SUC(G)$ see Lemmas
\ref{l:SUC-s-UC} and \ref{l:LE}.  
For the other assertions see \cite{GM}.
\end{proof}

\begin{remark}
The compactifications $G^{\mathrm{AP}}$ and $G^{\WAP}$ of $G$ are
respectively a {\em topological group} and a {\it semitopological
semigroup}. The compactifications $G^{\RUC}$ and $G^{\Asp}$ are
right topological semigroup compactifications of $G$ (see
\cite{GM}). The same is true for the compactification $j \colon G \to
G^{\SUC}$ (Lemma \ref{p:qLE}.5). Below (Theorem
\ref{t:conclusions}.5) we show that the {\it Roelcke
compactification} $i \colon G \hookrightarrow G^{UC}$ (which is always
proper by Lemma \ref{l:Roelcke}) is not in general a right
topological semigroup compactification. That is, $\UC(G)$ is not in
general point-universal.
\end{remark}

We sum up our results in the following dynamical hierarchy theorem
where we list dynamical properties of $f \in \RUC(X)$ and the
corresponding topological properties of $X_f \subset \RUC(G)$ (cf.
\cite[Remark~9.13]{GM}).

\begin{thm} \label{r:hierarchy} For every $G$-space $X$ and a
function $f \in \RUC(X)$ we have
\begin{align*}
&\text{$X_f$ is norm compact $\Longleftrightarrow$ $f$ is AP}\\
&\text{$X_f$ is weakly compact $\Longleftrightarrow$ $f$ is WAP}\\
&\text{$X_f$ is norm fragmented $\Longleftrightarrow$ $f$ is
Asplund
$\quad$ }\\
&\text{$X_f$ is orbitwise light $\Longleftrightarrow$ $f$
is LE $\quad$}\\
&\text{$X_f \subset UC(G)$  $\Longleftrightarrow$ $f$ is
SUC $\quad$}\\
\end{align*}
\end{thm}

\begin{example} \label{e:unit}
Let $G$ be the unitary group $U(H)=\Iso(H)$ where $H$ is an
infinite dimensional Hilbert space. Then
$UC(G)=\SUC(G)=\LE(G)=\Asp(G)=\WAP(G)$. Indeed the completion of $(G,
\mathcal L \bigwedge \mathcal R)$ can be identified with the
compact semitopological semigroup $\Theta(H)$ of all nonexpansive
linear operators (Uspenskij \cite{Us-un}). It follows that
$G^{\UC}$ can be identified with $\Theta(H)$.
The latter is a reflexively representable $G$-space (see for
example \cite[Fact 5.2]{Menz}). Therefore $\UC(G) \subset \WAP(G)$.
The reverse inclusion is well known (see for instance Theorem
\ref{l:inclusions}). Hence, $\UC(G) = \WAP(G)$.
\end{example}

\sk

Let $ \mathcal{C}_u=\{f\in UC(G): X_f\subset UC(G)\}.$ The
collection of functions $\mathcal{C}_u$ was studied by Veech in
\cite{V}. He notes there that $\WAP(G) \subset \mathcal{C}_u$ and
proves the following theorem.

\begin{thm} \label{t:veech}
\emph{(Veech \cite[Proposition 5.4]{V})} Let $G$ be a semisimple
analytic
Lie
 group with finite center and without compact factors. If $f\in
\mathcal{C}_u$ then every limit point of \ $Gf$ in $X_f$; i.e. any
function of the form $h(g)=\lim_{g_n\to\infty} f(gg_n)$, is a
constant function.
\end{thm}
\noindent (A sequence $g_n \in G$ ``tends to $\infty$" if each of
its ``projections" onto the simple components of $G$ tends to
$\infty$ in the usual sense.) He then deduces the fact that for
$G$ which is a direct product of simple groups the algebra
$\WAP(G)$ coincides with the algebra $\mathcal{W}^*$ of continuous
functions on $G$ which ``extend continuously" to the product of
the one-point compactification of the simple components of $G$
(\cite[Theorem 1.2]{V}). By our Proposition \ref{p:X_fin},
$\mathcal{C}_u = \SUC(G)$. Taking this equality into account,
Veech's theorem implies now the following result.

\begin{cor} \label{c:SL_n}
For every simple noncompact connected
Lie group $G$ with finite center (e.g.,
$\SL_n(\R)$) we have $\SUC(G)=\WAP(G)=\mathcal{W}^*$. 
In particular
the corresponding universal \emph{SUC} (and hence \emph{WAP}) compactification
is equivalent to the one point compactification of $G$.
\end{cor}

\section{The group $H_+[0,1]$}
\label{s:prop}

Consider the Polish topological group $G:=H_+[0,1]$ of all
orientation preserving homeomorphisms of the closed unit interval,
endowed with the compact open topology.
 Here is a list of some selected known results about this group:

\begin{enumerate} \label{r:H}
\item
$G$ is topologically simple.
\item
$G$ is not Weil-complete; that is, the right uniform structure
$\mathcal{R}$ of $G$ is not complete. The completion of the
uniform space $(G, \mathcal{R})$ can be identified with the
semigroup of all continuous, nondecreasing and surjective maps
$[0,1] \to [0,1]$ endowed with the uniform structure of uniform
convergence (Roelcke-Dierolf \cite[p. 191]{RD}).
\item
$G$ is \emph{Roelcke precompact} (that is the Roelcke uniformity
$\mathcal{L} \wedge \mathcal{R}$ on $G$ is precompact) \cite{RD}.
\item
The completion of $(G,\mathcal{L} \wedge \mathcal{R})$ can be
identified with the curves that connect the points $(0,0)$ and
$(1,1)$ and ``never go down"
(Uspenskij \cite{Uscurves}, see Lemma \ref{l:usp} below).
\item
Every weakly almost periodic function on $G$ is constant and every
continuous representation $G \to \Iso(V)$, where $V$ is a reflexive
Banach space, is trivial (Megrelishvili \cite{Merup}).
\item
$G$ is {\it extremely amenable}; that is every compact Hausdorff
$G$-space has a fixed point property
(Pestov \cite{pest-old}).
%
\end{enumerate}

We are going to show that $H_+[0,1]$ is SUC-trivial and hence also
${\mathcal Asp}$-trivial. Since every reflexive Banach space is
Asplund these results strengthen the main results of \cite{Merup}
(results mentioned in item (5) above).

\begin{defn}
Let $(X,\mu)$ be a compact $G$-space. We say that two points $a,b
\in X$ are SUC-{\it proximal} if there exist nets $s_i$ and $g_i$
in $G$ and a point $x_0 \in X$ such that $s_i$ converges to the
neutral element $e$ of $G$, the net $g_i x_0$ converges to $a$ and
the net $g_i s_i x_0$ converges to $b$.
\end{defn}

\begin{lem} \label{qle-proximal}
If the points $a$ and $b$ are SUC-proximal in a $G$-space $X$ then
$a \stackrel{\SUC} {\sim} b$.
\end{lem}
\begin{proof}
A straightforward consequence of our definitions using Lemma
\ref{l:SUC}.4.
\end{proof}

\begin{thm} \label{interval}
Let $G = H_+[0,1]$ be the topological group of
orientation-preserving homeomorphisms of $[0,1]$ endowed with the
compact open topology. Then $G$ is \emph{SUC}-trivial.
\end{thm}
\begin{proof}
Denote by $j \colon G \to G^{\SUC}$ and $i \colon G \to G^{UC}$ the
$G$-compactifications ($i$ necessarily is proper by Lemma
\ref{l:Roelcke}) induced by the Banach $G$-algebras $\SUC(G)
\subset \UC(G)$. There exists a canonical onto $G$-map $\pi \colon G^{UC}
\to G^{\SUC}$ such that the following diagram of $G$-maps is
commutative:

\begin{equation*}
\xymatrix { G \ar[dr]_{j} \ar[r]^{i} & G^{UC}
\ar[d]^{\pi} \\
  & G^{\SUC} }
\end{equation*}

We have to show that $G^{\SUC}$ is trivial for $G=H_+[0,1]$. One of
the main tools for the proof is the following identification.

\begin{lem}\label{l:usp}  \cite[Uspenskij]{Uscurves}
The dynamical system $G^{\UC}$ is isomorphic to the $G$-space
$(G,\Omega)$. Here $\Omega$ denotes the compact space of all
curves in $[0,1]\times [0,1]$ which connect the points $(0,0)$ and
$(1,1)$ and ``never go down", equipped with the Hausdorff metric.
These are the relations $\omega \subset [0,1] \times [0,1]$ where
for each $t\in [0,1]$, $\omega(t)$ is either a point or a vertical
closed segment. The natural action of $G=H_+[0,1]$ on $\Omega$ is
$(g\omega)(t) = g(\omega(t))$ (by composition of relations on
$[0,1]$).
\end{lem}


We first note that every ``zig-zag curve'' (i.e. a curve $z$ which
consists of a finite number of horizontal and vertical pieces) is
an element of $\Omega$. In particular the curves $\gamma_c $ with
exactly one vertical segment defined as $\gamma_c(t)=0$ for every
$t \in [0,c)$, $\gamma_c(c)=\{c\} \times [0,1]$ and $\gamma(t)=1$
for every $t \in (c,1]$,  are elements of $\Omega=G^{UC}$. Note
that the curve $\ga_1$ is a fixed point for the left $G$ action.
We let $\theta= \pi(\ga_1)$ be its image in $G^{\SUC}$. Of course
$\theta$ is a fixed point in $G^{\SUC}$. We will show that $\theta
= j(e)$ and since the $G$-orbit of $j(e)$ is dense in $G^{\SUC}$
this will show that $G^{\SUC}$ is a singleton.

The idea is to show that zig-zag curves are SUC-proximal in
$G^{\UC}$. Then Lemma \ref{qle-proximal} will ensure that their
images in $G^{\SUC}$ coincide. Choosing a sequence $z_n$ of zig-zag
curves which converges in the Hausdorff metric to $i(e)$ in
$G^{\UC}$ we will have $\pi(z_n)=\pi(\ga_1)=\theta$ for each $n$.
This will imply that indeed $j(e) = \pi(i(e)) = \pi(\lim_{n\to
\infty} z_n) = \lim_{n\to \infty} \pi(z_n)=\theta$.

First we show that $\pi(\gamma_1)=\pi(\gamma_c)$ for any $0<c<1$.
As indicated above, since $G^{\SUC}$ is the Gelfand space of the
algebra $\SUC(G)$, by Lemma \ref{qle-proximal}, it suffices to show
that the pair $\gamma_1, \gamma_c$ is SUC-proximal in $G^{UC}$.
Since $X^{\SUC}=X$ for $X:=G^{\SUC}$ we conclude that $\pi(\gamma_1)
= \pi (\gamma_c)$.

Let $p \in G^{\UC}$ be the
%
curve defined by $p(t)=t$ in the interval $[0,c]$ and by $p(t)=c$
for every $t \in [c,1)$. Pick a sequence $s_n$ of elements in $G$
such that $s_n$ converges to $e$ and $s_n c<c$. It is easy to
choose a sequence $g_n$ in $G$ such that $g_n s_n c$ converges to
$0$ and $g_nc$ converges to $1$. Then the sequences $s_n$ and
$g_n$ are the desired sequences; that is, $g_np \to \gamma_c$,
$g_n s_n p \to \gamma_1$ (see the picture below).
\sk

\centerline{\psfig{file=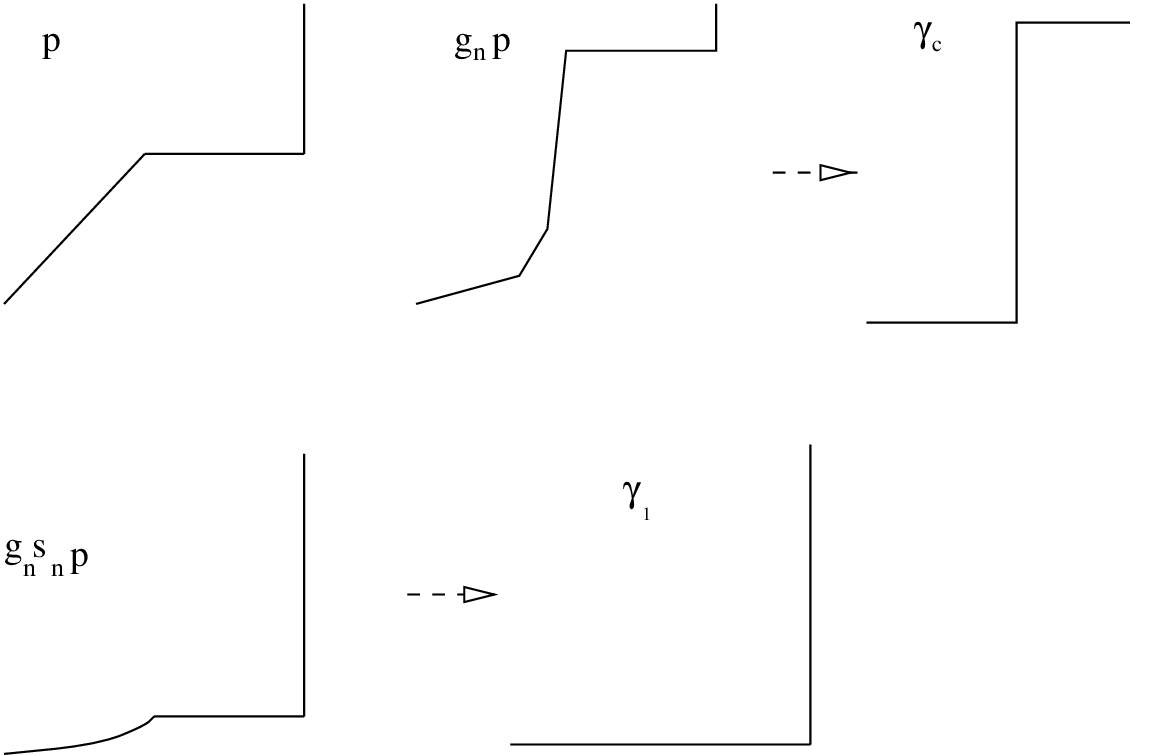,height=1.3in,width=8cm}}

\sk

Denote $\theta=\pi(\gamma_1)=\pi(\gamma_c)$. Using similar
arguments (see the picture below, where $a \overset{\pi}{\sim} b$
means $\pi(a)=\pi(b)$) construct a sequence $z_n \in G^{\UC}$ of
zig-zag curves which converges to $i(e)$ and such that $\pi
(z_n)=\theta$ for every $n$.

\sk

\centerline{\psfig{file=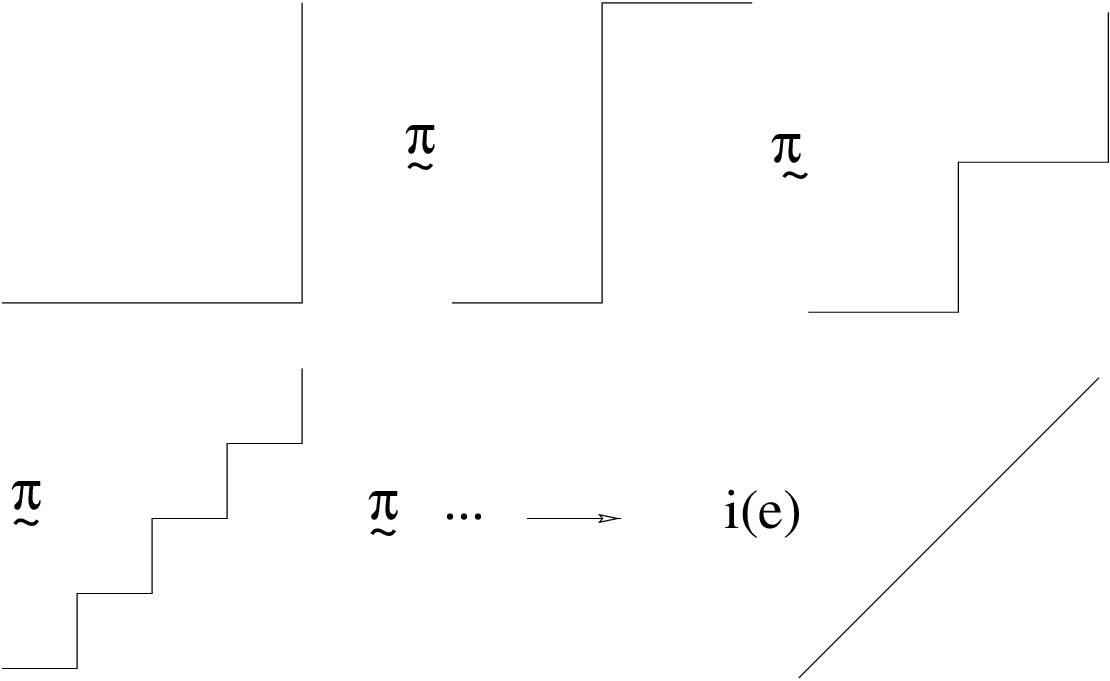,height=1.3in,width=8cm}}

\sk

In view of the discussion above this construction completes
the proof of the theorem.
\end{proof}

\section{Matrix coefficient characterization of SUC and LE}
\label{s:mat-SUC}

\begin{defn} \label{d:mat}
Let $h \colon G \to \Iso(V)$ be a co-representation of $G$ on a normed
space $V$ and let
$$
V \times G \to V, \ \ (v,g)
\mapsto vg:=h(g)(v)
$$
be the corresponding right action. For a pair of vectors $v \in V$
and $\psi \in V^*$ the associated \emph{matrix coefficient} is
defined by
$$
m_{v, \psi} \colon G \to \R, \hskip 0.3cm g \mapsto \psi(vg)=\langle vg,
\psi \rangle = \langle v, g \psi \rangle.
$$

If $h \colon G \to \Iso(V)$ is a representation then the matrix
coefficient $m_{v,\psi}$ is defined similarly by
$$
m_{v, \psi} \colon G \to \R, \hskip 0.3cm g \mapsto \psi(gv)=\langle gv,
\psi \rangle = \langle v, \psi g \rangle.
$$
\end{defn}

For example, if $V=H$ is a Hilbert space then $f=m_{u,\psi}$ is
the \emph{Fourier-Stieltjes transform}. In particular, for
$u=\psi$ we get the \emph{positive definite functions}.

We say that a vector $v \in V$ is \emph{$G$-continuous} if the
corresponding orbit map $\tilde{v} \colon G \to V, \ \tilde{v}(g)=vg$,
defined through $h \colon G \to \Iso(V)$, is norm continuous.
The \emph{continuous} $G$-vector $\psi \in V^*$
are defined similarly with respect to the dual action.

\begin{lem} \label{l:mat}
\cite{Menz} Let $h \colon G \to \Iso(V)$ be a co-representation of $G$ on
$V$. If $\psi$ (resp.: $v \in V$) is norm $G$-continuous, then $
m_{v,\psi}$ is left (resp.: right) uniformly continuous on $G$.
Hence, if $v$ and $\psi$ are both $G$-continuous then $ m_{v,\psi}
\in \UC(G)$.
\end{lem}

In the next theorem we list characterizations of several
subalgebras of $\RUC(G)$ in terms of matrix coefficients. These
characterizations also provide an alternative way to establish the
inclusions in Theorem \ref{l:inclusions}.2 for $X:=G$.

\begin{thm} \label{t:matrix}
Let $G$ be a topological group and $f \in C(G)$. \ben
\item
$f \in \RUC(G)$ iff $f=m_{v,\psi}$ for some continuous
co-representation $h \colon G \to \Iso(V)$, where $V \in Ban$.
\item
$f \in \UC(G)$ iff $f=m_{v,\psi}$ for some co-representation where
$v$ and $\psi$ are both $G$-continuous iff $f=m_{v,\psi}$ for some
continuous co-representation where $\psi$ is $G$-continuous.
\item
$f \in \SUC(G)$ iff $f =m_{v,\psi}$ for some continuous
co-representation $h \colon G \to \Iso(V)$, where $\varphi$ is norm
$G$-continuous in $V^*$ for every $\varphi$ from the weak star
closure $cl_{w^*}(G\psi)$.
\item
$f \in adj(G)$ iff $f=m_{v,\psi}$ for some adjoint continuous
co-representation.
\item $f \in \LE(G)$ iff $f =m_{v,\psi}$ for
some continuous co-representation $h \colon G \to \Iso(V)$, where weak
star and norm topologies coincide on each orbit $G\varphi$ where
$\varphi$ belongs to the weak star closure $Y:=cl_{w^*}(G\psi)$.
\item
$f \in \Asp(G)$ iff $f$ is a matrix coefficient of some continuous
Asplund co-representation of $G$.
\item
$f \in \WAP(G)$ iff $f$ is a matrix coefficient of some continuous
reflexive
co-representation (or, representation)
 of $G$.
 \een
\end{thm}

\sk

\begin{proof}
Claim (1) follows by taking in the regular
$\RUC(G)$-corepresentation $V:=\RUC(G)$, $v:=f$ and $\psi:=\a(e)$.
Claim (4) is a reformulation of Definition \ref{d:adj}.3.

The remaining assertions are essentially nontrivial. Their proofs
are based on an equivariant generalization of the
Davis-Figiel-Johnson-Pelczy\'nski interpolation technic
\cite{DFJP}. For detailed proofs of (6) and (7)
(for co-representations) see \cite[Theorem 7.17]{Menz} and
\cite[Theorem 5.1]{Menz}. As to the ``representations case" in (7)
observe that a matrix coefficient of
a (continuous) co-representation on a reflexive space $V$ can be
treated
as matrix coefficient of a (continuous) representation on the dual
space $V^*$. The continuity of the dual action follows by Theorem
\ref{t:adj}.

For (2) apply Lemma \ref{l:mat} and take $V:=\UC(G)$ (or, see  \cite{Meallmin}). Below we provide the proof of the new
assertions (3) and (5). See Theorems \ref{t:matSUC} and
\ref{t:LEmat} respectively.
\end{proof}

\begin{remark} \label{r:introverted}
Let $\Acal \subset \RUC(G)$ be a point-universal $G$-subalgebra.
Then $\Acal$ is \emph{left m-introverted} in the sense of
\cite{Mi}, \cite[Definition 1.4.11]{BJMo}. Indeed, we have only to
check that every matrix coefficient $m_{v,\psi}$ of the regular
$\Acal$-representation of the action $(G,G)$ on the Banach space
$\Acal$ belongs again to $\Acal$ whenever $v \in \Acal$ and $\psi
\in |\Acal| \subset \Acal^*$. Let $X:=|\Acal|$. Then the
$G$-system $(X,\eva(e))$ is point universal and $\Acal=\Acal(X,
\eva(e))$. The matrix coefficient $m_{v,\psi}$ comes from the
subsystem $({\cls}(G\psi), \psi)$. In other words $m_{v,\psi} \in
\Acal({\cls}(G\psi), \psi)$. By Lemma \ref{l:pu} we have
$\Acal({\cls}(G\psi), \psi) \subset \Acal(X, \eva(e))$. Thus
$m_{v,\psi} \in  \Acal(X, \eva(e))=\Acal$.
\end{remark}

\begin{defn} \label{d:small} \ 
\ben
\item
Let $G$ be a topological group and $G \times X \to X$ and $Y
\times G \to Y$ be respectively left and right actions. A map
$\langle \ ,\ \rangle \colon Y\times X \rightarrow \R$ is
$G$-\emph{compatible} if
$$
\langle y g,x\rangle=\langle y, gx\rangle \ \quad
\forall \ \ (y,g,x) \in Y \times G \times
X.
$$
\item
We say that a subset $M \subset Y$ is \emph{\emph{SUC}-small at} $x_0 \in
X$ if for every $\eps >0$ there exists a \nbd $U$ of $e$ such that
$$
 |\langle v, ux_0 \rangle -\langle v, x_0 \rangle |
\leq \eps \ \ \quad \forall \
(v,u) \in M \times U.$$
If $M$ is SUC-small at every $x \in X$
then we say that $M$ is \emph{\emph{SUC}-small for} $X$.
\item
Let $h \colon G \to \Iso(V)$ be a continuous co-representation on a
normed space $V$ and $h^* \colon G \to \Iso(V^*)$ be the dual
representation. Then we say that $M \subset V$ is \emph{\emph{SUC}-small
at} $x_0 \in X \subset V^*$ if this happens in the sense of (2)
regarding the canonical bilinear $G$-compatible map $\langle \ ,
\ \rangle \colon V \times V^* \to \R$. \een
\end{defn}

For example, a vector $\psi \in V^*$ in the dual space $V^*$ is
$G$-continuous iff the unit ball $B_V$ of $V$ is SUC-small at
$\psi$ (see Lemma \ref{l:small}.3).

We give some useful properties of SUC-smallness.

\begin{lem} \label{l:smallabs} \ 
\ben
\item Let $Y_1 \times X_1 \to \R$ and $Y_2 \times X_2 \to \R$
be two $G$-compatible maps. Suppose that $\g_1 \colon X_1 \to X_2$ and
$\g_2 \colon Y_2 \to Y_1$ are $G$-maps such that
$$
\langle y,\g_1(x)\rangle=\langle\g_2(y),x\rangle
\quad \quad \forall \ (y,x) \in Y_2
\times X_1.
$$
Then for every nonempty subset $M \subset Y_2$ the subset $\g_1(M)
\subset Y_1$ is \emph{SUC}-small at $x \in X_1$ if and only if $M$ is
\emph{SUC}-small at $\g_1(x) \in X_2$.
\item
Let $X$ be a (not necessarily compact) $G$-space. If $f \in
\SUC(X)$ then:
\bit
\item [(a)]
$fG$ is \emph{SUC}-small for $X$ \wrt the $G$-compatible evaluation map
$$fG \times X \to \R, \quad \ \ (fg,x) \mapsto f(gx).$$
\item [(b)]
The subset $F_e G$ of $C(X_f)$ is \emph{SUC}-small for $X_f$ considered
as a subset of $V^*$ where $V:=C(X_f)$ (with respect to the
canonical map $V \times V^* \to \R$
and the natural co-representation $G \to \Iso(V)$).
 \eit
 \een
\end{lem}
\begin{proof} (1): \
Observe that for every triple $(m, u, x_0) \in M \times G \times
X$ we have
$$
\langle \g_1(m), ux_0 \rangle - \langle \g_1(m),
x_0 \rangle = \langle m,\g_2(ux_0)\rangle - \langle m, \g_2(x_0)
\rangle=\langle m, u \g_2(x_0) \rangle - \langle m, \g_2(x_0)
\rangle.
$$

(2): \ (a) \ Directly follows by Lemma \ref{l:SUC}.4.

(b): \ Let $f \in \SUC(X)$. Then it comes by Proposition
\ref{f:X_f}.3 from a compact $G$-system $X_f$ and the
$G$-compactification $f_{\sharp} \colon X \to X_f$. As we know there
exists $F_e \in C(X_f)$ such that $f=F_e \circ f_{\sharp}$.
Theorem \ref{p:X_fin} implies that $X_f$ is SUC and $F_e \in
\SUC(X_f)$. By claim (a) it follows that $F_eG$ is SUC-small for
$X_f$. The $G$-compatible map $F_eG \times X_f \to \R$ can be
treated as a restriction of the canonical form $V \times V^* \to
\R$, where $V:=C(X_f)$ (considered $X_f$ as a subset of
$C(X_f)^*$).
\end{proof}

\begin{lem} \label{l:small}
Let $h \colon G \to \Iso(V)$ be a continuous co-representation.
\ben
\item
For every  subset $X$ of $V^*$ the family of \emph{SUC}-small sets for
$X$ in $V$ is closed under taking subsets, norm closures, finite
linear combinations,
finite unions
and convex hulls.
\item
If $M_n \subset V$ is \emph{SUC}-small at $x_0 \in V^*$ for every $n \in
\N$ then so is the set $\bigcap_{n \in \N} (M_n + \delta_n B_V)$
for every positive decreasing sequence $\delta_n \to 0$.
\item  For every $\psi \in V^*$ the following are
equivalent:
\bit
\item [(i)]
The orbit map $\widetilde{\psi} \colon G \to V^*$ is norm continuous.
\item [(ii)]
$\mathcal B$ is SUC-small at $\psi$, where 
$$\mathcal
B:=\{\breve{v} \colon V^* \to \R, \ x \mapsto \check{v}(x):=\langle v,x
\rangle \}_{v \in B_V}.$$ 
\eit
\een
\end{lem}
\begin{proof}
Assertion (1) is straightforward.

(2): \ We have to show that the set $\bigcap_{n \in \N} (M_n +
\delta_n B_V)$ is SUC-small at $x_0$. Let $\eps
>0$ be fixed. Since $Gx_0$ is a bounded subset of $V^*$ one can
choose $n_0 \in \N$ such that $|v(gx_0)| < \frac{\eps}{4}$ for
every $g \in G$ and every $v \in \delta_{n_0} B_V$. Since
$M_{n_0}$ is SUC-small at $x_0$ we can choose a \nbd $U(e)$ such
that $|m(ux_0)-m(x_0)| < \frac{\eps}{2}$ for every $u \in U$ and
every $m \in M_{n_0}$. Now every element $w \in \bigcap_{n \in \N}
(M_n + \delta_n B_V)$ has a form $w=m + v$ for some $m \in
M_{n_0}$ and $v \in \delta_{n_0} B_V$. Then for every $u \in U$ we
have
$$|w(ux_0) - w(x_0)| \leq |m(ux_0)- m(x_0)| + |v(ux_0)| + |v(x_0)|
< \frac{\eps}{2} + \frac{\eps}{4}+ \frac{\eps}{4} = \eps.
$$

(3): \ Use that $\norm{u\psi-\psi}=sup_{v \in B_V} |\langle v,
u\psi \rangle -\langle v, \psi\rangle|$ and $B_V$ is
$G$-invariant.
\end{proof}

\begin{thm} \label{t:matSUC}
The following conditions are equivalent:
\ben
\item
$f \in \SUC(G)$.
\item
$f =m_{v,\psi}$ for some continuous Banach co-representation $h \colon G
\to \Iso(V)$,
$v\in V$ and $\psi \in V^*$,
with the property that
$\varphi$ is norm $G$-continuous
for every $\varphi$ in the weak $^*$ closure $cl_{w^*}(G\psi)$.
 \een

Moreover, one can assume in \emph{(2)} that $cl_{w^*}(G\psi)$
separates points of $V$.
\end{thm}

\begin{proof} (2) $\Longrightarrow$ (1):
Let $h \colon G \to \Iso(V)$ be a continuous co-homomorphism such that
$f$ is a matrix coefficient of $h$. That is, we can choose $v \in
V$ and $\psi \in V^*$ such that $f(g)=\langle
vg,\psi\rangle=\langle v,g \psi \rangle$ for every $g \in G$. One
can assume that $||\psi||=1$. The strong continuity of $h$ ensures
that the dual restricted (left) action of $G$ on the weak star
compact unit ball $(V^*_1,w^*)$ is jointly continuous. Consider
the orbit closure $X:={\cls}_{w^*} (G\psi)$ in the compact
$G$-space $(V^*_1,w^*)$. Define the continuous function ${\hat v}:
X \to \R$ induced by the vector $v$. Precisely, ${\hat
v}(x)=\langle v,x \rangle$ and in particular, $f(g)={\hat
v}(g\psi)$. So $f$ comes from $X$ and the compactification $\nu \colon G
\to X, \ g \mapsto g\psi$. It suffices to show that the $G$-system
$X$ is SUC. Let $x_0$ be an arbitrary point in $X$ and let $w$ be
an arbitrary vector in $V$. By the definition of the weak star
topology and the corresponding uniformity on the compact space $X$
it suffices to show that for every $\eps >0$ there exists a \nbd
$U(e)$ such that $|{\hat w} (gux_0) -{\hat w}(gx_0)| \leq \eps$
for every $g \in G$ and every $u \in U$. By simple computations we
get
\begin{align*}
|{\hat w} (gux_0) -{\hat w}(gx_0)| &=|\langle w,gux_0\rangle
-\langle w,gx_0 \rangle| =
|\langle wg,ux_0 \rangle - \langle wg,x_0 \rangle|\\
&=|\langle wg,ux_0-x_0 \rangle| \leq ||wg|| \cdot ||ux_0 - x_0||.
\end{align*}

Take into account that $||wg||=||w||$. Since $x_0 \in
{\cls}(G\psi)$, by our assumption the orbit map $\widetilde{x_0}:
G \to V^*$ is norm continuous with respect to the dual action of
$G$ on $V^*$. Therefore, given $\varepsilon > 0$, there exists a
neighborhood $U$ of $e$ in $G$ such that $||ux_0 - x_0|| <
||w||^{-1}\eps$ for every $ u \in U$. Thus, $|{\hat w} (gux_0)
-{\hat w}(gx_0)| \leq \eps$. This shows that $X$ is SUC. Hence, $f
\in \SUC(G)$.

(1) $\Longrightarrow$ (2): Let $f \in \SUC(X)$. Then it comes by
Proposition \ref{f:X_f}.3 from a compact transitive $G$-system
$X_f$ and the $G$-compactification $f_{\sharp} \colon G \to X_f$. There
exists $F:=F_e \in C(X_f)$ such that $f=F \circ f_{\sharp}$.
By Lemma \ref{l:smallabs}.2 we conclude that $FG$ is SUC-small for
$X_f \subset C(X_f)^*$.

Let $M:=co(-FG \cup FG)$ be the convex hull of the symmetric set \
$-FG \cup FG$. Then $M$ is a convex symmetric bounded
$G$-invariant subset in $C(X_f)$. By Lemma \ref{l:small}.1 we know
that $M$ is also SUC-small for $X_f$.

For brevity of notation let $E$ denote the Banach space $C(X_f)$.
Since $X_f$ is a compact $G$-space the natural right action of $G$
on $E=C(X_f)$ (by linear isometries) is continuous.

Consider the sequence $K_n:=2^n M + 2^{-n} B_E$, where $B_E$ is
the unit ball of $E$. Since $M$ is convex and symmetric, we can
apply the construction of \cite{DFJP} (we mostly use the
presentation and the development given by Fabian in the book
\cite{Fa}). Let $\| \ \|_n$ be the Minkowski's functional of the
set $K_n$. That is,
$$
\| v\|_n = \inf\ \{\lambda
> 0 \bigm| v\in \lambda K_n\}
$$
Then $\| \ \|_n$ is a norm on $E$ equivalent to the given norm
of $E$ for every $n \in \N$. For $v\in E,$ let

$$
N(v):=\left(\sum^\infty_{n=1} \| v \|^2_n\right)^{1/2}
\hskip 0.1cm \text{and}
\hskip 0.1cm \hskip 0.1cm V: = \{ v \in E \bigm| N(v) < \infty
\}
$$

Denote by $j \colon V \hookrightarrow E$ the inclusion map. Then $(V,
N)$ is a Banach space, $j \colon V \to E$ is a continuous linear
injection and
$$
M \subset j(B_V)=B_V
$$
Indeed, if
$v \in M$
 then $2^nv \in K_n$. Therefore, $\| v\|_n \leq 2^{-n}$ and
 $N(v)^2 \leq \sum_{n \in \N} 2^{-2n} <1$.

\vskip 0.1cm

By our construction $M$ and $B_E$ are $G$-invariant. This implies
that the natural right action $V \times G \to V, \ \ (v,g) \mapsto
vg$ is isometric, that is $N(vg)=N(v)$. Moreover, by the
definition of the norm $N$ on $V$ (use the fact that the norm
$\norm{\cdot}_n$ on $E$ is equivalent to the given norm of $E$ for
each $n \in \N$) we can show that this action is norm continuous.
Therefore, the co-representation $h \colon G \to \Iso(V), \ h(g)(v):=vg$
on the Banach space $(V,N)$ is well defined and continuous.

Let $j^* \colon E^* \to V^*$ be the adjoint map of $j \colon V \to E$.
Now our aim is to check the $G$-continuity of every vector
$\varphi \in j^*(X_f)=cl_{w^*}(G\psi)$, where $\psi:=j^*(z)$ and
$z$ denotes the point $f_{\sharp}(e) \in X_f$.
By Lemma \ref{l:small}.3 we have to show that $B_V$ is SUC-small
for $j^*(X_f)$.

\sk

\nt \textbf{Claim} : $j(B_V) \subset \bigcap_{n \in \N} K_n =
\bigcap_{n \in \N} (2^n M + 2^{-n}B_E)$.
\begin{proof}
The norms $\norm{\cdot}_n$ on $E$ are equivalent to each other. It
follows that if $v \in B_V$ then $\| v\|_n < 1$ for all $n \in
\N$. That is, $v \in \lambda_n K_n$ for some $0 < \lambda_n <1$
and $n \in \N$. By the construction $K_n$ is a convex subset
containing the origin. This implies that $\lambda_n K_n \subset
K_n$. Hence $j(v)=v \in K_n$ for every $n \in\N$.
\end{proof}


Recall now that $FG$ is SUC-small for $X_f \subset C(X_f)^*$. By
Lemma \ref{l:small}.1 we know that also $M:=co(-FG \cup FG)$ is
SUC-small for $X_f \subset C(X_f)^*$. Moreover by Lemma
\ref{l:small}.2 we obtain that $A:=\bigcap_{n \in \N} (2^n M +
2^{-n}B_E) \subset C(X_f)$ is SUC-small for $X_f \subset
C(X_f)^*$. The linear continuous operator $j \colon V \to C(X_f)$ is a
$G$-map. Then by Lemma \ref{l:smallabs}.1 it follows that
$j^{-1}(A) \subset V$ is SUC-small for $j^*(X_f) \subset V^*$. The
same is true for $B_V$ because by the above claim we have $j(B_V)
\subset A$ (and hence, $B_V \subset j^{-1}(A)$). That is $B_V$ is
SUC-small for $j^*(X_f)$. Now Lemma \ref{l:small}.3 shows that the
orbit map $\widetilde{\varphi} \colon G \to V^*$ is $G$-continuous for
every
$\varphi \in j^*(X_f)=cl_{w^*}(G\psi)$. By our construction $F \in
j(V)$ (because $F \in M \subset j(B_V) $). Since $j$ is injective
the element $v:=j^{-1}(F)$ is uniquely determined in $V$.
We show that $f=m_{v,\psi}$ for the co-representation $h$. Using
the equality $F \circ \a_f =f$ and the fact that $\a_f$ is a
$G$-map we get
$$
\langle Fg, z \rangle = F(g\a_f(e))=(F \circ \a_f)(g)=f(g).
$$
On the other hand,
$$
m_{v,\psi}(g)=\langle vg,\psi \rangle = \langle j^{-1}(F)g, j^*(z)
\rangle = \langle j(j^{-1}(F))g),z \rangle = \langle Fg, z
\rangle.
$$
Hence, $f=m_{v,\psi}$, as required. Therefore we have proved that
(1) $\Longleftrightarrow$ (2).

Finally we show that one can assume in (2) that
$cl_{w^*}(G\psi)=j^*(X_f)$ separates points of $V$. If $v_1, v_2$
are different elements in $V$ then $j(v_1) \neq j(v_2)$. Since
$X_f$ separates $C(X_f)$ then $\langle j(v_1), \phi \rangle \neq
\langle j(v_2), \phi \rangle$ for some $\phi \in X_f$. Now observe
that $\langle j(v),\phi \rangle = \langle v, j^*(\phi) \rangle$
for every $v \in V$.
\end{proof}

\begin{cor} \label{c:adjisSUC}
${\mathcal Adj}(G) \subset \SUC(G)$.
\end{cor}

\sk

Next we show how one can characterize $\LE(G)$ in terms of matrix
coefficients.

\begin{thm} \label{t:LEmat}
The following conditions are equivalent:
\ben
\item
$f \in \LE(G)$.
\item
$f =m_{v,\psi}$ for some continuous co-representation $h \colon G \to
\Iso(V)$,
$v\in V$ and $\psi \in V^*$, for a Banach space $V$,
with the property that the weak $^*$ and the norm
topologies coincide on the orbit $G\varphi$ of every $\varphi$ in
the weak $^*$ closure $Y:=cl_{w^*}(G\psi)$.
\een

Moreover, one can assume in (2) that $Y$ separates points of $V$.
\end{thm}
\begin{proof}
(2) $\Longrightarrow$ (1): \ By definition $f$ comes from
$Y:=cl_{w^*}(G\psi)$. Hence it suffices to show that $Y$ is LE.
Equivalently, we need to show that $Y$ is orbitwise light (see
\cite[Lemma 5.8.2]{GM}). Let $\mu_Y$ be the uniform structure on
the compact space $Y$. Denote by $(\mu_Y)_G$ the corresponding
uniform structure of uniform convergence inherited from $Y^G$ (see
\cite{GM}). We have to show that
$top(\mu_Y)|_{G\varphi}=top((\mu_Y)_G)|_{G\varphi}$ for every
$\varphi \in Y$. Observe that the topology $(\mu_Y)_G$ on the
orbit $G\varphi$ is weaker than the norm topology. Since the
latter is the same as the weak star topology (that is,
$top(\mu_Y)|_{G\varphi}$) we get that indeed
$top(\mu_Y)|_{G\varphi}= top((\mu_Y)_G)|_{G\varphi}$.

(1) $\Longrightarrow$ (2): \
The proof uses again the interpolation technique of \cite{DFJP},
as in Theorem \ref{t:matSUC}. The proof is similar so we omit the
details. However we provide necessary definition and two lemmas
(Definition \ref{dd:small} and Lemmas \ref{l:LEsmallnessabs} and
\ref{l:LEsmallness}). They play the role of Lemmas
\ref{l:smallabs} and \ref{l:small}.

For every set $M$ denote by $\R^M$ the set of all real valued
functions $M \to \R$. The
topologies of pointwise
and uniform convergence on $\R^M$ will be denoted by $\tau_p$ and
$\tau_u$ respectively.

\begin{defn}\label{dd:small}
Let $\langle \ , \ \rangle \colon Y\times X \rightarrow \R$ be a
$G$-compatible map (as in Definition \ref{d:small}) and $M$ be a
nonempty subset of $Y$. Denote by $j \colon X \to \R^M, \
j(x)(m):=\langle m,x\rangle$ the associated map.
 \ben
\item
We say that a subset $A$ of $X$ is \emph{$M$-light} if the
pointwise and uniform topologies coincide on
$j(A) \subset \R^M$.
\item
 $M$ is
\emph{LE-small at} $x_0 \in X$ if the orbit $Gx_0$ is $M$-light.
\item
$M$ is \emph{LE-small for} $X$ if the orbit $Gx$ is $M$-light
at every $x \in X$ (compare Theorem \ref{l:LEand}). \een
\end{defn}

 We are going to examine this definition in a
particular case of the canonical bilinear map $V \times V^* \to
\R$ which is $G$-compatible for every co-representation $h \colon G \to
\Iso(V)$.

We collect here some useful properties of LE-smallness.

\begin{lem} \label{l:LEsmallnessabs} \ 
\ben
\item Let $Y_1 \times X_1 \to \R$ and $Y_2 \times X_2 \to \R$
be two $G$-compatible maps. Suppose that $\g_1 \colon X_1 \to X_2$ and
$\g_2 \colon Y_2 \to Y_1$ are $G$-maps such that
$$
\langle y,\g_1(x)\rangle=\langle
\g_2(y),x\rangle  \quad \quad \forall \ (y,x) \in Y_2
\times X_1.
$$
Then for every nonempty subset $M \subset Y_2$ the subset $\g_1(M)
\subset Y_1$ is LE-small at $x \in X_1$ if and only if $M$ is
LE-small at $\g_1(x) \in X_2$.
\item
Let $X$ be a (not necessarily compact) $G$-space. If $f \in \LE(X)$
then the subset $F_e G$ of $C(X_f)$ is LE-small for $X_f$
considered as a subset of $V^*$ where $V:=C(X_f)$ (with respect to
the canonical map $V \times V^* \to \R$
and the natural co-representation $G \to \Iso(V)$).
 \een
\end{lem}
\begin{proof}
(1): \ Similar to Lemma \ref{l:smallabs}.1.

(2): \ Let $f \in LE(X)$. Then it comes by Proposition
\ref{f:X_f}.3 from a compact $G$-system $X_f$ and the
$G$-compactification $f_{\sharp} \colon X \to X_f$. As we know there
exists $F_e \in C(X_f)$ such that $f=F_e \circ f_{\sharp}$.
Theorem \ref{l:LEand} implies that $X_f$ is LE and $F_e \in
LE(X_f)$. By the same theorem, $X_f$ is orbitwise light in
$\RUC(G)$. This means that pointwise and norm topologies in
$\RUC(G)$ agree on every $G$-orbit in $X_f$. 
On the other hand it is straightforward to see that for the
$G$-compatible map
$$
F_eG \times X_f \to \R
$$
(Definition \ref{dd:small} with $M:=F_eG$) the corresponding
pointwise topology $\tau_p$ on
$X_f$ coincides with the pointwise topology inherited from
$\RUC(G)$ and the uniform topology $\tau_u$ on $X_f$ coincides with
the norm topology of $\RUC(G)$.
\end{proof}

\begin{lem} \label{l:LEsmallness}
Let $h \colon G \to \Iso(V)$ be a continuous co-representation.
 \bit
\item [(a)]
For every $X \subset V^*$ the family
of LE-small sets for $X$ in $V$ is closed under taking:
subsets, norm closures, finite linear combinations,
finite unions 
and convex hulls.
\item [(b)]
If $M_n \subset V$ is LE-small at $x_0 \in V^*$ for every $n \in
\N$ then so is 
$$\bigcap_{n \in \N} (M_n + \delta_n B_V)$$
 for every
positive decreasing sequence $\delta_n \to 0$.
\item [(c)]
The following are equivalent:
\bit
\item [(i)]
The pointwise and norm topologies agree on the $G$-orbit $Gx$ for
every $x \in X \subset V^*$.
\item [(ii)] $\mathcal B$ is LE-small for
$X$, where $$\mathcal B:=\{\breve{v} \colon V^* \to \R, \ x \mapsto
\check{v}(x):= \langle v,x \rangle \}_{v \in B_V}.$$ 
\eit
 \eit
\end{lem}
\end{proof}

\begin{question} 
Do we have ${\adj}(G)=\SUC(G)$ for every topological group $G$
? The question seems to be open even for abelian
non-discrete $G$ (say $G=\R$). 
\nl The equivalent
question for abelian $G$ is whether ${\adj}(G)=\UC(G)$ ?
Also, what is the relation between the algebras
$\LE(G)$ and $\adj(G)$ ?
\end{question}


\section{Some conclusions about $H_+[0,1]$ and $\Iso(\U_1)$}
\label{s:cor}

From the reflexive triviality of $H_+[0,1]$ and results of
Uspenskij about $\Iso(\U_1)$ Pestov deduces in \cite[Corollary
1.4]{Pe-new} the fact that the group $\Iso(\U_1)$ is also
reflexively trivial. Using a similar idea and the matrix
coefficient characterization of SUC and LE one can conclude that
$\Iso(\U_1)$ is SUC-trivial and LE-trivial.

Recall the following results of Uspenskij.

\begin{thm} \label{f:usp} \emph{(Uspenskij
\cite{Us-subgr,UspComp})}
 The group $\Iso(\U_1)$ is topologically simple and contains a
copy of every second countable topological group (e.g.,
$H_+[0,1]$).
\end{thm}

Lemma \ref{l:observation} below is a generalized version of
Pestov's observation.
Of course it is important here that the
corresponding property admits a reformulation in terms of Banach
space representations,
which
is the case, for instance, for SUC and LE.

\begin{lem} \label{l:observation}
Let $G_1$ be a topological subgroup of a group $G_2$. Suppose that
$G_2$ is $G_1$-\emph{simple}, in the sense that, every
\emph{non-trivial} normal subgroup $N$ in $G_2$ containing $G_1$
is necessarily dense in $G_2$. Then if $G_1$ is either: 1)
SUC-trivial, 2) LE-trivial, 3) adjoint continuous trivial or 4)
$\mathcal{K}$-trivial (where $\mathcal{K}$ is a class of Banach
spaces) then the same is true for $G_2$.
\end{lem}
\begin{proof}
We consider only the case of SUC. Other cases are similar
(and even easier for (3) and (4)).

We use Theorem \ref{t:matSUC}. Let $h \colon G_2 \to \Iso(V)$ be a
continuous co-representation where $Y:=cl_{w^*}(G\psi)$ separates
points of $V$ and $\varphi$ is norm $G_2$-continuous in $V^*$ for
every $\varphi \in Y$. It is enough to show that any such
co-representation of $G_2$ is trivial. By Theorem \ref{t:matSUC}
this will show that $G_2$ is SUC-trivial. First observe that the
restriction $h|_{G_1}$ of $h$ to $G_1$ is trivial. In fact,
otherwise $vg \neq v$ for some $(v,g) \in V \times G_1$ and by our
assumption there exists $\varphi \in Y$ such that $\varphi(v) \neq
\varphi(vg)$. Then the restriction $m_{v,\varphi}|_{G_1}$ of the
corresponding matrix coefficient $m_{v,\varphi} \colon G_2 \to \R$ to
$G_1$ is not constant. However, by Theorem \ref{t:matSUC},
$m_{v,\varphi}|_{G_1} \in \SUC(G_1)$, contradicting our assumption
that $G_1$ is SUC-trivial. Therefore, $h|_{G_1} \colon G_1 \to \Iso(V)$
is trivial. Hence $G_1$ is a subgroup of the normal closed
subgroup $N:=ker(h)$ of $G_2$. Since $G_2$ is $G_1$-simple it
follows that $N=G_2$. Hence $h$ is trivial.
\end{proof}

Note that if $G_2$ is \emph{topologically simple}, that is
$\{e\}$-simple, then it is $G_1$-simple for every subgroup.

The following theorem sums up some of our results concerning the
topological groups $H_+[0,1]$ and $\Iso(\U_1)$.

\begin{thm} \label{t:conclusions}
Let $G$ be one of the groups $H_+[0,1]$ or $\Iso(\U_1)$. \ben
\item
The compactifications
$G^{\SUC}, G^{\LE}, G^{\Asp}, G^{\WAP}$ are trivial.
\item
Every adjoint continuous (co)representation of the group $G$ is
trivial.
\item
Every continuous Asplund (co)representation of the group $G$ is
trivial.
\item Every (co)representation $h \colon G \to \Iso(V)$
on a \emph{separable} Asplund space $V$ is trivial.
\item
The algebra $\UC(G)$ and the ambit $(G^{\UC},i(e))$ are not point
universal. In particular, the map $i \colon G \to G^{UC}$ is not a right
topological compactification of $G$. \een
\end{thm}
\begin{proof} (1): $H_+[0,1]$ is SUC-trivial
by Theorem \ref{interval}. By results of Uspenskij (see Theorem
\ref{f:usp}) the group $\Iso(\U_1)$ is topologically simple and
also a universal second countable group. In particular it contains
a copy $G_1$ of $H_+[0,1]$ as a topological subgroup. It follows
that $G_2:=\Iso(\U_1)$ is $G_1$-simple. Applying Lemma
\ref{l:observation} we conclude that $G_2:=\Iso(\U_1)$ is also
SUC-trivial. The rest follows by the inclusions of Theorem
\ref{l:inclusions}.

(2): Every adjoint continuous (co)representation of $H_+[0,1]$
must be trivial. Otherwise, by Theorem \ref{t:matSUC} (or
Corollary \ref{c:adjisSUC}), it contains a nonconstant SUC
function. Now Lemma \ref{l:observation} implies that $\Iso(\U_1)$
is also adjoint continuous trivial.

(3): By Theorem \ref{t:adj} every continuous Asplund
(co)representation of $G$ is adjoint continuous. Now apply (2).

(4): By a recent result of Rosendal and Solecki \cite[Corollary
3]{RS} every homomorphism of $G=H_+[0,1]$ into a separable group
is necessarily continuous. Combining this result and our assertion
(3) we obtain a proof in the case of $G=H_+[0,1]$. The case of
$G=\Iso(\U_1)$ now follows by using again the $H_+[0,1]$-simplicity
of $\Iso(\U_1)$.

(5): Take a non-constant uniformly continuous function on $G$
(such a function necessarily exists by Lemma \ref{l:Roelcke}).
Since $\SUC(G)=\{\text{constants}\}$ we get $\SUC(G) \neq \UC(G)$.
Now Corollary \ref{p:CUSUC} finishes the proof.
\end{proof}


By Theorems \ref{f:usp} and \ref{t:conclusions} we get

\begin{cor}
Every second countable group $G_1$ is a subgroup of a Polish
\emph{SUC}-trivial group $G_2$.
\end{cor}

However the following questions are open (see also \cite{opit}).

\begin{question} \ 
\ben
\item
Find a nontrivial 
Polish group which is \emph{SUC}-trivial (${\mathcal Ref}$-trivial,
${\mathcal Asp}$-trivial) but does not contain a subgroup
topologically isomorphic to $H_+[0,1]$.
\item
Is the group $H(I^{\omega})$ \emph{SUC}-trivial (${\mathcal
Ref}$-trivial, ${\mathcal Asp}$-trivial)?

\emph{And, a closely related question
(see Lemma \ref{l:observation}):}
\item
Is the group $G_2:=H(I^{\omega})$, $G_1$-simple for a subgroup
$G_1 < G_2$ where $G_1$ is a copy of either $H_+[0,1]$ or of
$\Iso(\U_1)$? \een
\end{question}

\begin{thm} \label{t:noMETRICcomp} Let $G$ be an Asplund trivial
(e.g. $H_+[0,1]$ or $\Iso(\U_1)$) group. Then every
\emph{metrizable} right topological semigroup compactification of
$G$ is trivial.
\end{thm}
\begin{proof}
By Theorem \ref{t:conclusions}.1, $G^{\Asp}$ is trivial, so that
every RN transitive $G$-space is trivial. If $G \to S$ is a right
topological semigroup compactification of $G$, then the natural
induced $G$-space $(G,S)$ is isomorphic to its own enveloping
semigroup. By a recent work \cite{GMU}, a metric dynamical system
$(G,X)$ is RN iff its enveloping semigroup is metrizable. Now if
$S$ is metrizable then it follows that the transitive system
$(G,S)$ is RN and therefore trivial. 
\end{proof}

%
%
%
%

Recall, in contrast, that for every topological group $G$ the
algebra $\RUC(G)$ separates points and closed subsets on $G$ and
therefore the maximal right topological semigroup compactification
$G \hookrightarrow G^{\RUC}$ is
faithful.

\section{Relative extreme amenability: SUC-fpp groups} \label{s:amenable}

Recall that a topological group $G$ has the {\em fixed point on
compacta property (fpp)} (or is \emph{extremely amenable}) if
every compact 
$G$-space $X$ has a fixed point. It is
well known that locally compact extremely amenable groups are
necessarily trivial
(see for example \cite{GL}). 
Gromov and Milman \cite{GrMi} proved that the unitary group $U(H)$
is extremely amenable. Pestov has shown that the groups $H_+[0,1]$
and $\Iso(\U_1)$ are extremely amenable (see \cite{pest-old,
Pe-book} for more information).

Consider the following relativization.

\begin{defn} \label{d:fpp}
Let P be a class of compact $G$-spaces. \ben
\item
A $G$-space $X$ is P-{\it fpp}
(or is {\em extreme P-amenable})
if every $G$-compactification $X
\to Y$ such that $Y$ is a member of P has a fixed point.
\item
A topological group $G$ is P-{\it fpp}
(or is {\em extremely P-amenable})
if the $G$-space $X:=G$ is
P-fpp or equivalently, if every $G$-space $Y$ in P has a fixed
point.
 \een
\end{defn}

Taking P as the collection of all compact flows we get extreme
amenability. With the class P of compact affine flows we recover
amenability. When P is taken to be the collection of
equicontinuous (that is, almost periodic) flows we obtain the old
notion of minimal almost periodicity (MAP). Minimal almost
periodicity was first studied by von Neumann and Wigner \cite{vNW}
who showed that $\mathrm{PSL}(2,\Q)$ has this property. See also Mitchell
\cite{Mi} and Berglund, Junghenn and Milnes \cite{BJMo}.

\begin{lem}
Let P be a class of compact $G$-spaces which is preserved by
isomorphisms, products subsystems and quotients. Let ${\mathcal
P}$ and $X^{{\mathcal P}}$ be as in Section \ref{s:actions}. The
following conditions are equivalent: \ben
\item
$G$ is P-fpp.
\item
The compact $G$-space $G^{{\mathcal P}}$ is P-fpp.
\item
Any minimal compact $G$-space in P is trivial.
\item
For every $f \in {\mathcal P}$ the $G$-system $X_f$ has a fixed
point.
\item
The algebra ${\mathcal P}$ is extremely left amenable (that is it
admits a multiplicative left invariant mean). \een
\end{lem}
\begin{proof}
Clearly each of the conditions (1) and (3) implies all the others.
Use the fact that the $G$-space $G^{{\mathcal P}}$ is point
universal to deduce that each of (2) and (5) implies (3). Finally
(4) implies (2) because $G^{{\mathcal P}}$ has a presentation as a
subsystem of the product of all the $X_f$, $f\in \mathcal{P}$.
Note that $f \in {\mathcal P}$ iff $X_f$ has property $P$ (see
\cite[Proposition 2.9.3]{GM}).
\end{proof}

\begin{remark}\label{class} \ 
\ben
\item
The smaller the class P
is one expects the property of being P-fpp to be less restrictive;
however even when one takes P to be the class of equicontinuous
$\Z$-spaces (that is, \emph{cascades}) it is still an open
question whether P-fpp, that is minimal almost periodicity, is
equivalent to extreme amenability
(see \cite{G}).
\item
A minimal compact $G$-space $X$ is LE iff $X$ is AP. It follows by
the inclusions $\rm{LE} \supset \rm{RN_{app}} \supset WAP \supset
AP$ (cf. also Theorem \ref{l:inclusions}) that $G$ is minimally
almost periodic iff $G$ is P-fpp for each of the following
classes: WAP, $\rm{RN_{app}}$ or LE.
\een
\end{remark}

Here we point out two examples of topological groups $G$ which are
SUC-extremely amenable (equivalently, SUC-fpp)
but not extremely amenable. In the next
two sections we will show
that $S_{\infty}$ as well as the group $H(C)$ of homeomorphisms of
the Cantor set $C$ are also SUC-fpp (both groups are not extremely
amenable).
See Corollaries \ref{env-cor} and   
\ref{c:H(C)} below.

\begin{example} \label{e:lc}
For every $n \geq 2$ the simple Lie group $\SL_n(\R)$, being
locally compact, is not extremely amenable. However it
is SUC-extremely amenable.
%
This follows easily from Corollary \ref{c:SL_n}.
\end{example} 

\section{The Roelcke compactification of the group
$S(\N)$}\label{Sec-S}

Let $G=S(\N)$ be the Polish topological group of all permutations
of the set $\N$ of natural numbers (equipped with the topology of
pointwise convergence). Consider the one point compactification
$X^*=\N \cup \{\infty\}$ and the associated natural $G$ action
$(G,X^*)$. For any subset $A\subset\N$ and an injection $\alpha:A
\to \N$ let $p_{\alpha}$ be the map in $(X^*)^{X^*}$ defined by
$$
p_\alpha(x) =
\begin{cases}
\alpha(x) & \quad x\in A\\
\infty & \quad{\text{otherwise}}
\end{cases}
$$
We have the following simple  claim.
\begin{claim}
The enveloping semigroup $E=E(G,X^*)$ of the $G$-system $(G,X^*)$
consists of the maps $\{p_\alpha \colon A \to \Z\}$ as above.
Every element of $E$ is a continuous function so that by the
Grothendieck-Ellis-Nerurkar theorem \cite{EN}, the system
$(G,X^*)$ is \emph{WAP}.
\end{claim}

\begin{proof}
Let $\pi_\nu$ be a net of elements of $S(\N)$ with
$p=\lim_{\nu}\pi_\nu$ in $E$. Let $A=\{n\in \N \colon p(n)\ne\infty\}$
and $\alpha(n)=p(n)$ for $n\in A$. Clearly $\alpha \colon A \to \N$ is
an injection and $p=p_\alpha$.

Conversely given $A\subset \N$ and an injection $\alpha:A \to\N$
we construct a sequence $\pi_n$ of elements of $S(\N)$ as follows.
Let $A_n=A \cap [1,n]$ and $M_n=\max\{\alpha(i) \colon i\in A_n\}$. Next
define an injection $\beta_n:[1,n] \to \N$ by
$$
\beta_n(j) =
\begin{cases}
\alpha(j) & \quad j\in A\\
j+ M_n +n & \quad{\text{otherwise}}.
\end{cases}
$$
Extending the injection $\beta_n$ to a permutation $\pi_n$ of
$\N$, in an arbitrary way, we now observe that
$p_{\a}=\lim_{n\to\infty}\pi_n$ in $E$. The last assertion is
easily verified.
\end{proof}

\begin{thm} \label{env} \ 
\begin{enumerate}
\item
The two algebras $\UC(G)$ and $\WAP(G)$ coincide for $G=S(\N)$. 
\item
The universal \emph{WAP} 
compactification $G^{\WAP}$ of $G$ (and hence also $G^{\UC}$), is
isomorphic to $E=E(G,X^*)$.
Thus the universal \emph{WAP} (and
Roelcke) compactification of $G$ is homeomorphic to the Cantor
set.
\end{enumerate}
\end{thm}  
\begin{proof}
Given $f\in \UC(G)$ and an $\eps >0$ there exists $k\in \N$ such
that --- with $H=H(1,\dots,k)=\{g\in G: g(j)=j,\ \forall \ 1 \le j
\le k\}$ ---
$$
\sup_{u, v \in H}|f(ugv) - f(g)| < \eps.
$$
Set $\hat f(g) = \sup_{u, v \in H}f(ugv)$, then
$\|\hat f - f\| \le \eps$.
Clearly $\hat f$, being $H$-biinvariant, is both
right and left uniformly continuous; i.e. $\hat f\in \UC(G)$.
Let
$$
\N^k_*=\{(n_1,n_2,\dots,n_k) : n_j\in \N \text{\ are distinct}\}
=\{\text{injections}: \ \{1,2, \cdots, k\} \to \N \}
$$
and let $G$ act on $\N^k_*$ by
$
g(n_1,n_2,\dots,n_k)=(g^{-1}n_1,g^{-1}n_2,\dots,g^{-1}n_k).
$
The stability group of the point $(1,\dots,k)\in \N^k_*$ is just
$H$ and we can identify
the
discrete $G$-space $G/H$ with $\N^k_*$.
Under this identification, to a function $f\in \UC(G)$ which is
right $H$-invariant (that is $f(gh)=f(g),\ \forall g\in G, h\in
H$), corresponds a bounded function $\om_f \in
\Omega_k=\R^{\N^k_*}$, namely
$$
\om_f(n_1,n_2,\dots,n_k) = f(g) \quad \text{iff} \quad \
g(j)=n_j,\ \forall\ 1 \le  j \le k.
$$
If we now assume that $f\in \UC(G)$ is both right and left
$H$-invariant (so that $f=\hat f$) then, as we will see below,
$f$ and accordingly its corresponding $\om_f$, admits only
finitely many values, corresponding to the
finitely many double $H$ cosets $\{HgH: g\in G\}$.

We set $Y_f=Y={\cls}\{g\om_f: g\in G\}\subset
\Omega_k=\R^{\N^k_*}$, where the closure is with respect to the
pointwise convergence topology. $(G,Y_f)$ is a compact $G$-system
which is isomorphic, via the identification $G/H\cong \N^k_*$, to
$X_f\subset \R^G$.
We will refer to elements of $\Omega_k=\R^{\N^k_*}$ as {\em
configurations}. Consider first the case $k=2$.

In the following figure we have a representation of the
configuration $f=\om_f=\om_{1,2}$ and three other typical elements
of $Y_f$. The configuration $\om_{2,7}= \sigma\om_{1,2}$ --- where
$\sigma$ is the permutation $\left(
\begin{smallmatrix}
1 & 2 & 7\\
7 & 1 & 2
\end{smallmatrix}
\right) $
--- , admits seven values (the maximal number it
can possibly have): ``blank" at points $(m,n)$ with $m,n \not\in
\{2,7\}$, the values $\diamond$ and $*$ at $(2,7)$ and $(7,2)$
respectively, and four more constant values on the two horizontal
and two vertical lines. (The circled diagonal points $(2,2)$ and
$(7,7)$ are by definition not in $\N^2_*$.) If we let $\pi_n$ be
the permutation
$$
\pi_n(j) =
\begin{cases}
j & \quad j \not\in \{1,n\}\\
n & \quad {\text{for}}\ j=1\\
1 & \quad {\text{for}}\ j=n
\end{cases}
$$
and denote by $p=\lim_{n\to\infty}\pi_n$, the corresponding
element of $E(G,X)$
then, e.g.
$\om_{1,\infty}=p\om_{1,2}=\lim_{n\to\infty}\pi_n\om_{1,2}=
\lim_{n\to\infty}\om_{1,n}$.

\br


\centerline{\psfig{file=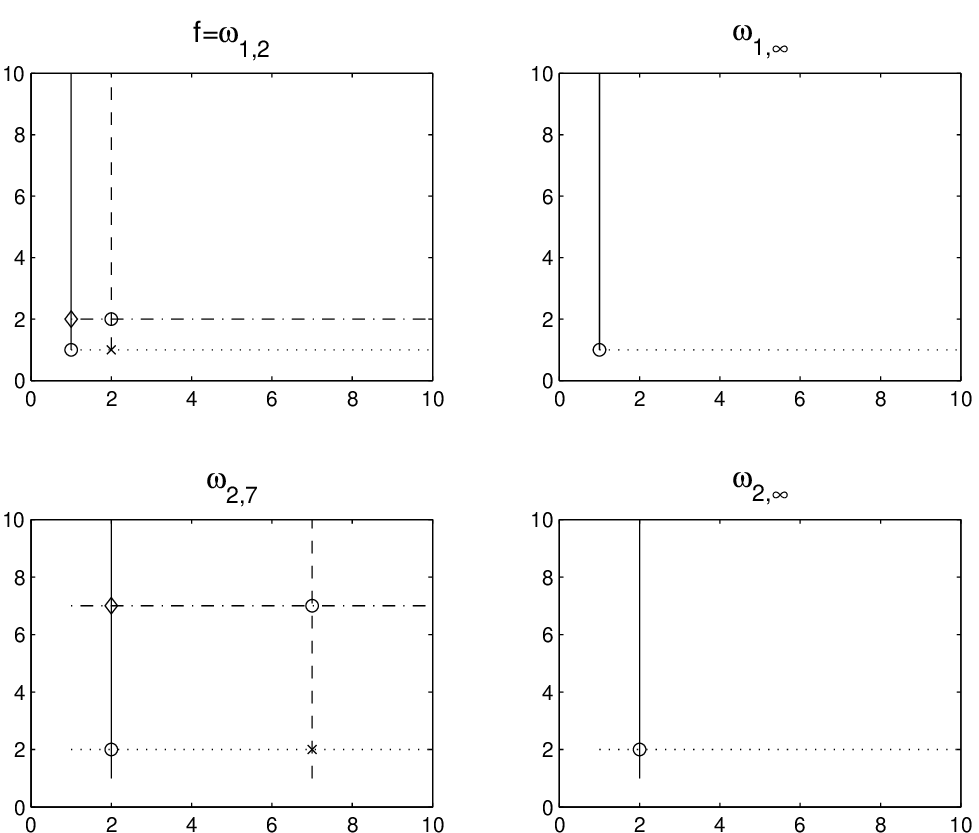,height=2in,width=6cm}}
\centerline{The functions $\omega_{1,2}$, $\omega_{1,\infty}$,
	$\omega_{2,7}$ and $\omega_{2,\infty}$ in $\Omega_2$}

\br



Now it is not hard to see that the $G$ action on $Y$ naturally
extends to an action of $E=E(G,X^*)$ on $Y$ where each $p\in E$
acts continuously.
(Show that the map $(a,b) \mapsto \om_{a,b}$ is an isomorphism of
$G$-systems from $X^*\times X^*\setminus \Del$ onto $Y$, where
$\Del=\{(n,n): n\in \N\}$.) It then follows that $E=E(G,X^*)$
coincides with $E(G,Y)$.

By the Grothendieck-Ellis theorem (see e.g. \cite[Theorem
1.45]{G1}) these observations show that the $G$-space $Y_f$ is WAP
and therefore the function $f$, which comes from $(G,Y_f)$, is a
WAP function.

These considerations are easily seen to hold for any positive
integer $k$. For example, an easy calculation shows that for
$H_k=H(1,2,\dots,k)$ the number of double cosets $\{H_kgH_k: g\in
G\}$  is
$$
\sum_{j=0}^k \binom{k}{j}\frac{k!}{(k-j)!}
$$

Since the subgroups $H_k=H(1,\dots,k)$ form a basis
for the topology at $e$, as we have already seen, the union of the
$H_k$-biinvariant functions for $k=2,3,\dots$ is dense in $\UC(G)$
and we conclude that indeed $\UC(G)=\WAP(G)$.

Since for each $H$-biinvariant function $f$ the enveloping
semigroup of the dynamical system $(G,Y_f)$ is isomorphic to
$E(G,X^*)$ and since the Gelfand compactification of $UC(G)$ is
isomorphic to a subsystem of the direct product
$$
\prod \{Y_f: f \ \text{is $H$-biinvariant for some
$H(1,\dots,k)$}\}
$$
we deduce that $E(G,X^*)$ serves also as the enveloping semigroup
of the universal dynamical system $|\WAP(G)|=|\UC(G)|$. Finally,
since $|\WAP(G)|$ is point universal we conclude (by Lemma
\ref{l:pu}) that $(G,|\WAP(G)|)$ and $(G,E(G,X^*))$ are
$G$-isomorphic.
\end{proof}

\begin{cor}\label{env-cor}
The Polish group $G=S(\N)$ 
is \emph{SUC}-extremely amenable 
but not extremely amenable.
\end{cor} 
\begin{proof}
It was shown by Pestov  \cite{pest-old} that $G$ is not extremely
amenable and the nontrivial universal minimal $G$-system was
described in \cite{GW3}. On the other hand the $G$-system $G^{\UC}$
described in Theorem \ref{env} admits a unique minimal
set which is a fixed point. Thus
the SUC $G$-system $G^{\SUC}$, being
$G$-isomorphic to $G^{\UC}$ (see Theorems \ref{env} and
\ref{l:inclusions}.2), has a fixed point.
\end{proof}

\section{The homeomorphisms group of the Cantor set}

In this section let $C$ denote the classical Cantor set --- i.e.
the ternary subset of the interval $[0,1]$. Thus $C$ has the
following representation:
$$
C=\bigcap_{n=0}^\infty I^n,
$$
where $I^n=\bigcup_{j=1}^{2^n} I^n_j$, is the disjoint union of
the $2^n$ closed intervals obtained by removing from $I=[0,1]$ the
appropriate $2^n -1$ open `middle third' intervals. In the sequel
we will write $I^m_j$ for the clopen subset $I^m_j\cap C$ of $C$.
For each integer $m\ge 1$, $\mathbf{I}^m=\{I^m_j: 1\le j \le
2^m\}$ denotes the basic partition of $C$ into $2^m$ clopen
``intervals".

We let $G=H(C)$ be the Polish group of homeomorphisms of $C$
equipped with the topology of uniform convergence. For $n\in\N$ we
let
$$
H_n=\{g\in G: g I^n_j =I^n_j,\ \forall\  1\le j \le 2^n\}.
$$
Each $H_n$ is a clopen subgroup of $G$ and we note that the system
of clopen subgroups $\{H_n: n=2,3,\dots\}$ forms a basis for the
topology of $G$ at the identity $e\in G$.

For any fixed integer $k\ge 1$ consider the collection
$$
\Acal^k=\{\mathbf{a}=\{A_1,A_2,\dots,A_k\}: {\text{ a partition of
$C$ into $k$ nonempty clopen sets}}\}.
$$
In particular note that for $k=2^n$, $\mathbf{I}^n$ is an element
of $\Acal^k$.

The discrete homogeneous space $G/H_n$ can be identified with
$\Acal^{k}=\Acal^{2^n}$: an element $gH_n  \in  G/H_n$ is uniquely
determined by the partition
$$
\mathbf{a}=\{gI^n_j:  1\le j \le 2^n = k\},
$$
and conversely to every partition $\mathbf{a}\in \Acal^{k}$
corresponds a coset $gH_n \in G/H_n$. In fact, if
$\mathbf{a}=\{A_1,A_2,\dots,A_k\}$ we can choose $g$ to be any
homeomorphism of $C$ with $A_j=gI^n_j$.

Thus for $k=2^n$ we have a parametrization of $\Acal^k$ by the
discrete homogeneous space $G/H_n$.

Let
$$
\Om^k=\R^{\Acal^k} \cong \R^{G/H_n}.
$$
Via the quotient map $G \to G/H_n,\ g\mapsto gH_n$, the Banach
space $\ell^\infty(\Acal^k)$ canonically embeds into the Banach
space $\RUC(G)$ where the image consists of all the right
$H_n$-invariant functions in $\RUC(G)$. Thus if $f\in \RUC(G)$
satisfies $f(gh)=f(g)$ for all $g\in G$ and $h\in H_n$ then
$\om_f(gI^n_1,\dots,gI^n_k)=\om_f(A_1,\dots,A_k)=f(g)$, where
$A_j=gI^n_j$, is the corresponding configuration in
$\ell^\infty(\Acal^k)$.

We equip $\Om^k=\R^{\Acal^k}$ with its product topology. The group
$G$ acts on the space $\Om^k$ as follows. For $\om\in \Om^k$ and
$g\in G$ let
$$
g\om(\mathbf{a})=\om(g^{-1}A_1,g^{-1}A_2,\dots,g^{-1}A_k),
$$
for any $\mathbf{a}=\{A_1,A_2,\dots,A_k\}\in \Acal^k$.
Equivalently $g\om(g'H_n)=\om(g^{-1}g'H_n)$, for every $g'H_n\in
G/H_n$. For each right $H_n$-invariant functions $f$ in $\RUC(G)$
we denote the compact orbit closure of $f=\om_f$ in $\Om^k$ by
$Y_f$.

First let us consider the case $n=1$, where $k=2$,
$$
\Acal^2=\{\mathbf{a}=\{A,A^c\}: {\text{ a partition of $C$ into
$2$ nonempty clopen sets}}\},
$$
and
$$
H=H_1=\{g\in G: g I^1_j =I^1_j,\ \  j=1,2\}.
$$  
\begin{claim}
There are exactly seven double cosets $HgH$, $g\in G$.
\end{claim}  
\begin{proof}
For a partition $(A,A^c)\in \Acal^2$ exactly one of the following
five possibilities holds: (1) $A=I^1_1$, (2) $A=I^1_2$, (3) $A
\varsubsetneqq I^1_1$, (4) $ A \subsetneqq I^1_2$, (5) $A
\supsetneqq I^1_1$, (6) $A \supsetneqq I^1_2$, (7) $A\cap I^1_1
\ne \emptyset \ne A\cap I^1_2$, and $A^c\cap I^1_1 \ne \emptyset
\ne A^c\cap I^1_2$.

Clearly for any two partitions $(A,A^c), (B,B^c)$ we have
$(B,B^c)=(hA,hA^c)$ for some $h\in H$ iff they belong to the same
class. Our claim follows in view of the correspondence $G/H \cong
\Acal^2$.
\end{proof}

Define an element $\om_f\in \Om^2$ and the corresponding function
$f\in \UC(G)$ as follows:
$$
\om(A,A^c)=j \qquad \text{if $(A,A^c)$ is of type ($j$)},
\qquad j=1,\dots,7,
$$
and $f(g)=\om(gI^1_1,gI^1_2)$. Clearly $f$ is $H_1$-biinvariant
and in particular an element of $UC(G)$. Let $X_f$ denote the
(pointwise) orbit closure of $f$ in $\RUC(G)$. Via the natural lift
of $\Om^2$ to $\RUC(G)$ we can identify $X_f$ with $Y_f = {\cls}
\{g\cdot\om_f: g\in G\}\subset \Om^2$.

Next consider a sequence of homeomorphisms $h_n\in G$ satisfying
the conditions
\begin{enumerate}
\item[(i)]
$h_n(I^n_1)=(I^n_{2^n})^c$,
\item[(ii)]
$h_n((I^n_1)^c)=I^n_{2^n}$ and
\item[(iii)]
$h_n$ is order preserving.
\end{enumerate}
It is then easy to check that the limit $\lim_{n\to\infty}h_n\om_f
= \om_0$ exists in $\Om^2$ where $\om_0$ is defined by
$$
\om_0(A,A^c)=
\begin{cases}
5 & \text{if $0 \in A$}\\
4 & \text{if $0 \not\in A$}.
\end{cases}
$$
Now for any $g\in G$ we have
$$
(g\cdot\om_0)(A,A^c)= \om_0(g^{-1}A,g^{-1}A^c)=
\begin{cases}
5 & \text{if $g(0) \in A$}\\
4 & \text{if $g(0) \not\in A$}.
\end{cases}
$$
For $x\in C$ set
$$
(\om_x)(A,A^c)=
\begin{cases}
5 & \text{if $x \in A$}\\
4 & \text{if $x \not\in A$}.
\end{cases}
$$
Then for $g\in G$ we have $g\om_0=\om_{g0}$. Moreover denoting
$Y_0 = {\cls} \{g \om_0: g\in G\}\subset \Om^2$ we have
$Y_0=G\om_0$ and the map $\phi \colon (G,C) \to (G,Y_0)$ defined by
$\phi(x)=\om_x$ is an isomorphism of $G$-spaces. We get the
following lemma.

\begin{lem}\label{finite}
Let $Y_0={\cls}\{g \om_0: G\in G\}$ be the orbit closure of
$\om_0$ in $\Om^2$, then the $G$-space $(G,Y_0)$ is isomorphic to
$(G,C)$, the natural action of $G=H(C)$ on the Cantor set $C$.
\end{lem}

\begin{remark}
An argument analogous to that of Lemma \ref{finite} will show that
for every $n$ the number of $H_n$ double cosets is finite. As in
the case of $S(\mathbb{N})$ in the previous section this shows the
well known fact  that $G=H(C)$ is Roelcke precompact
(see \cite{Uscurves}).
\end{remark}

\br

In contrast to Theorem \ref{env} we obtain the following result.

\begin{thm} \label{t:H(C)}
For $G=H(C)$ we have $\UC(G)\supsetneqq \SUC(G).$
\end{thm}
\begin{proof}



Consider the function
$$
f_0(g)=\om_0(g^{-1}I^1_1,g^{-1}I^1_2)=g\omega_0(I^1_1,I^1_2)
=\omega_{g0}(I^1_1,I^1_2)
$$
and let $h_n\in G$ be defined as above. Let $u_n$ be a sequence of
elements of $G$ which converges to $e\in G$ and for which $h_n u_n
0 = 2/3$. Then, as $h_n 0 = 0$ for every $n$, we have
$$
f_0(h_n) = \om_0(h_n^{-1}I^1_1,h_n^{-1}I^1_2)=5,
$$
but as $h_n u_n 0 = 2/3$,
$$
f_0(h_n u_n)  = \om_0(u_n^{-1}h_n^{-1}I^1_1,u_n^{-1}h_n^{-1}I^1_2)
= \om_{h_n u_n 0}(I^1_1,I^1_2) = 4.
$$
Thus $f_0$ is not left uniformly continuous. Since $f_0 \in X_f
\cong Y_0$, we conclude, by Theorem \ref{p:X_fin}, that $f$ is not
a SUC function.
\end{proof}

\begin{remark} \label{cantorSUC}
A similar argument will show that any two points $a,b\in C$ are
SUC-proximal for the $G$-space $(G,C)$. Thus this $G$-space is
SUC-trivial by Lemma \ref{qle-proximal}. Letting $F:Y_0 \to
\{4,5\}\subset \R$ be the evaluation function
$F(\om)=\om(I^1_1,I^1_2)$, we observe that
$$
f_0(g)=\om_0(g^{-1}I^1_1,g^{-1}I^1_2)=g\omega_0(I^1_1,I^1_2)
=F(g\om_0) =
F (g \phi 0) = (F\circ\phi)(g 0).
$$
Thus the function $f_0$ comes from the $G$ space $C$, via the
continuous function $F\circ \phi \colon C \to \R$ and the point $0\in
C$. This is another way of showing that $f_0$ and hence also $f$
are not SUC.
\end{remark}


\begin{remark} 
By Theorem \ref{t:H(C)} and Corollary \ref{p:CUSUC} we obtain, in
particular, that the algebra $\UC(G)$ is not point universal and
the corresponding Roelcke compactification $G \to G^{\UC}$ is not a
right topological semigroup compactification of $G$. The same is
true for $G:=H_+[0,1]$ because $\UC(G) \neq \SUC(G)$.
This follows from Theorem
\ref{interval} and Lemma \ref{l:Roelcke}.

\end{remark}

\section{Topological 2-transitivity 
	vs SUC} 


%
	%
%

\begin{defn}
	Let $\pi \colon G \times X \to X$ be a continuous action. 
	
	\begin{enumerate}
		
		\item  A point $x \in X$ is 
		\textit{transitive} if the $G$-orbit $Gx = \{gx: g \in G\}$ is dense in $X$. We denote by 
		$\mathrm{Trans}(X)$ the set of transitive points in $X$.  
		The action is called {\em point transitive} when $ \mathrm{Trans}(X)$ is nonempty.
		\item   The action is called {\em topologically transitive} when the set
		$\{g \in G: gO_1 \cap O_2\}$ is nonempty for every pair of nonempty open subsets $O_1, O_2$ in $X$.
		\item
		We say that the action on $X$ is {\em weakly mixing} when the diagonal action on
		$X^2 = X \times X$ is topologically transitive.	 	
		
		\item 	A point $(x_1,x_2) \in X^2$ is 
		2-\textit{transitive} if the $G$-orbit $\{(gx_1,gx_2): g \in G\}$ of $(x_1,x_2)$ is dense in $X^2$. We denote by 
		$\mathrm{Trans_2}(X)$ the set of $2$-transitive points in $X^2$.  
		
		\item The action $\pi$ is \emph{algebraically $2$-transitive} if the induced diagonal action 
		$$G \times X^2 \setminus \Delta \to X^2 \setminus \Delta$$ is algebraically transitive 
		(one orbit).

	\end{enumerate} 
\end{defn}


When $X$ is compact and metrizable and the action is topologically transitive, 
then the subset $\mathrm{Trans}(X)$ of all transitive points is an invariant 
dense $G_{\delta}$-subset of $X$ 
(indeed, $\mathrm{Trans}(X)=\cap_{n \in \N} (\cup \{g^{-1}O_n: g \in G\})$ 
for every countable base $\{O_n: n \in \N\}$ of $X$). 
Thus, in this case,  $X$ is topologically transitive if and only if $X$ is point transitive. 
Therefore, if $X$ is a compact metrizable $G$-flow, then $X$ is weakly mixing
if and only if $ \mathrm{Trans_2}(X)$ is nonempty.
Note that if $\mathrm{Trans_2}(X)$ is nonempty then it is a $G$-invariant dense subset of $X^2$. 
If in addition $X$ has an isolated point, then $X$ is a singleton.

\begin{prop} \label{p:top2trans} 
	Let $\pi \colon G \times X \to X$ be a continuous action. Assume that there exists $c_0 \in X$ such that:
	\begin{enumerate}
		\item its $G$-orbit $Gc_0$ is dense in $X$;
		\item the orbit map $G \to Gc_0$ is open;
		
		\item $(Gc_0 \times Gc_0) \cap \mathrm{Trans_2}(X)$ is dense in $X \times X$.
	\end{enumerate}
	Then the $G$-space $X$ is \emph{SUC}-trivial (hence, $G$ has \emph{SUC}-fpp).   
\end{prop}
\begin{proof} 
	It is enough to show that every $f \in \SUC(X)$ is constant. 
	Assume to the contrary that $f \colon X \to \R$ is a
	nonconstant SUC function.
	Then there exist: $\eps >0$ and $x_1,x_2 \in X$ such that 
	$$
	|f(x_1)-f(x_2)|>\eps. 
	$$
	By the density assertion (3) and the continuity of $f$, there exists $(w_1,w_2) \in (Gc_0 \times Gc_0) \cap \mathrm{Trans_2}(X)$ which is sufficiently close to $(x_1,x_2)$ in $X^2$ such that 
	$$
	|f(w_1)-f(w_2)|>\eps. 
	$$
	By Lemma \ref{l:SUC}.4 
	there exists an open neighborhood $U(e)$ in $G$ such that $U^{-1}=U$ and 
	$$
	|f(guc_0)-f(gc_0)| < \frac{\eps}{2} \ \ \ \forall (g,u) \in G \times U.
	$$
	The triangle inequality ensures that  
	\begin{equation} \label{eq:inMix} 
		|f(gu_1c_0)-f(gu_2c_0)| < \eps \ \ \ \forall (g,u_1,u_2) \in G \times U \times U.	
	\end{equation} 
	By condition (2) the image $Uc_0$ is open in $Gc_0$. Let $O$ be an open subset of $X$ such that $O \cap Gc_0 =Uc_0$. 
	
	Since $(w_1,w_2) \in \mathrm{Trans_2}(X)$, 
	we can choose $g_0 \in G$ such that 
	$$(g_0w_1,g_0w_2) \in O \times O.$$
	Then, since $(w_1,w_2) \in Gc_0 \times Gc_0$, we have  
	$$(g_0w_1,g_0w_2) \in (O \cap Gc_0) \times (O \cap Gc_0)= Uc_0 \times Uc_0.$$
	Therefore, there exist $u_1, u_2 \in U$ such that $w_1=g_0^{-1}u_1c_0$ and $w_2=g_0^{-1}u_2c_0$ for some $u_1, u_2 \in U$. Now, Equation \ref{eq:inMix} implies that 
	$$
	|f(w_1)-f(w_2)|<\eps, 
	$$
	which contradicts
	the choice of $(w_1,w_2)$. 
\end{proof}

\begin{thm} \label{t:GeneralCor} 
	Let $G$ be a Polish topological group such that $M(G)$ is metrizable and proximal. Then $M(G)$ is \emph{SUC}-trivial 
	(hence, $G$ has the \emph{SUC}-fpp). 	
\end{thm}
\begin{proof}
	Every minimal proximal flow is weakly mixing by 
	\cite[Chapter 2, Cor. 2.2]{Gl-book}. Hence, $M(G)$ is weakly mixing. 
	As we already mentioned, 
	when $M(G)$ is compact metric, it is weakly mixing if and only if 
	there exists a point of transitivity for the diagonal action of $G$ on $M(G) \times M(G)$. 
	Then, $\mathrm{Trans_2}(M(G))$ is a dense $G_{\delta}$-subset of $M(G) \times M(G)$. 
	
	On the other hand, since the universal minimal $G$-flow $M(G)$ is metrizable, by a criterion due to  
	Ben Yaacov,  Melleray and Tsankov \cite{BMT} 
	(which for nonarchimedean groups $G$ previously was proved by Zucker \cite{Zucker14}), 
	there exists a \textit{generic point} $c_0 \in M(G)$. That is, we have a dense $G$-orbit $Gc_0 \subseteq M(G)$ which is a $G_{\delta}$-subset of $M(G)$. By Effros' theorem, the orbit map $G \to Gc_0, g \mapsto gc_0$ is open. 
	Then, the set $Gc_0 \times Gc_0$ and also $(Gc_0 \times Gc_0) \cap \mathrm{Trans_2}(M(G))$ are dense $G_{\delta}$-subsets of $M(G) \times M(G)$. Now Proposition \ref{p:top2trans} finishes the proof. 	
\end{proof}

\begin{cor} \label{c:H(C)} 
	Let $G=H(C)$ be the Polish homeomophism group of the Cantor set  $C$. Then $M(G)$ is \emph{SUC}-trivial 
	(hence, $G$ has the \emph{SUC}-fpp). 	 
\end{cor}
\begin{proof}
	By   \cite{GW4},  
	$M(H(C))$ is metrizable and proximal. 
\end{proof}

\begin{remark} 
	Corollary \ref{c:H(C)} was formulated in 
	\cite[Theorem 13.8]{Gl-Me-SUC}. However, its proof was not correct. 
	We thank Lionel Van Th\'{e} for pointing out this error.  
\end{remark}

\begin{cor} \label{c:Aut} 
	Let $G:=\Aut(\Q,\circ)$ be 
	the Polish group of all circular order preserving permutations of $\Q_0$ with the pointwise topology, where  $\Q_0$ is the rational discrete circle. Then $M(G)$ is \emph{SUC}-trivial. 	
\end{cor}
\begin{proof}
	By \cite[Thm 5.2]{GM-UltraHom}, $M(G)$ is a circularly ordered metrizable compactum which we get from the circle after splitting its rational points. 
	Also, $M(G)$ is a proximal $G$-flow by \cite[Thm 4.9]{GM-UltraHom}.  	
\end{proof}

\sk  
If the action $G \times X \to X$ is algebraically $2$-transitive and $X$ is perfect 
(i.e., has no isolated points) then the action is 
topologicaly 2-transitive and hence 
Proposition \ref{p:top2trans} applies. 
For many concrete homogeneous metric compact perfect 
spaces $X$ the natural action of the topological group $G=H(X)$ on
$X$ is algebraically 2-transitive. By Proposition \ref{p:top2trans} the flow
$(G,X)$ admits only constant SUC functions and the corresponding
SUC $G$-compactification $X^{\SUC}$ is trivial. This is the case,
to mention some concrete examples, for $X$ the Cantor set, the Hilbert cube and the circle $\mathbb{T}$.

In the latter case even the subgroup 
$G:=H_+(\mathbb{T}) < H(\mathbb{T})$ of all
orientation-preserving homeomorphisms of the circle acts
algebraically 2-transitive on $\mathbb{T}$.
Pestov has shown \cite{pest-old, Pe-book} 
that the universal minimal dynamical $G$-system $M(G)$ for
$G:=H_+(\mathbb{T})$ coincides with the natural action of $G$ on
$\mathbb{T}$.

Combining these results with Proposition \ref{p:top2trans} we obtain 
\begin{cor}
	The Polish group $G=H_+(\mathbb{T})$ of orientation preserving
	homeomorphisms of the circle has the \emph{SUC}-fpp (is \emph{SUC}-extremely amenable)   
	but it is not extremely amenable.
\end{cor}

This can be proved also using Corollary \ref{c:Aut} because there exists a (dense) injective continuous homomorphism $\Aut(\Q,\circ) \to H_+(\mathbb{T})$.

An alternative proof follows easily from Proposition \ref{hom}.
Moreover, the following result is stronger than Proposition \ref{hom} and leads to an additional explanation of Theorem \ref{t:GeneralCor}. 

\begin{thm} \label{t:stronger}
	Let $X$ be a compact minimal $G$-space which contains a topologically transitive point $x_0 \in X$ such that the orbit map $G \to Gx_0$ is open. 
	\ben
	\item
	If $X$ is \emph{SUC} then $X$ is AE. 
	\item If $X$ is \emph{SUC} and minimal then $X$ is 
	equicontinuous. 
	\een
\end{thm}
\begin{proof} (1): 
	Since the orbit map $G \to gx_0$ is open, the orbit $Gx_0$ can be identified with the coset $G$-space $G/H$, with the stabilizer subgroup $H=st(x_0)$.  
	As in the proof of Proposition \ref{hom} we can verify that $x_0$ is a point of equicontinuity of the action of $G$ on $(G/H,\mu)$, where $\mu$ is the $\mu$ precompact uniformity on $G/H$ inherited from the compact space $X$. Then 
	$x_0$ is a point of equicontinuity also in $cl(Gx_0)=X$ (Lemma \ref{l:LE}.1). 
	So, $Eq(X)$ is nonempty. In particular, $X$ is non-sensitive. 
	By Lemma \cite[Lemma 9.2.3]{GM}, $\mathrm{Trans}(X) \subseteq Eq(X)$. 
	Clearly, $\mathrm{Trans}(X)$ is dense in $X$ (containing $Gx_0$). 
	Therefore, $Eq(X)$ is dense in $X$. This means that $X$ is AE.
	
	(2) Follows from (1) because $\mathrm{Trans}(X) \subseteq Eq(X)$. By the minimality, $X=\mathrm{Trans}(X)$. Hence, $X=Eq(X)$. 	 
\end{proof}

Note again that 
 for Polish groups $G$ with metrizable $M(G)$ there exists a generic point $x_0 \in M(G)$. For every $G$-factor $q \colon M(G) \to Y$ the point $q(x_0)$ is generic in $Y$ (see \cite[Prop. 14.1]{AKL}). These results, together with Theorem  \ref{t:stronger} imply that every SUC $G$-factor of a proximal metric $M(G)$ is trivial (being proximal and equicontinuous). This gives one more verification of Theorem \ref{t:GeneralCor}. 


\begin{remark} 
	The	SUC-fpp for topological groups can be 
	otherwise described as \emph{extremal SUC-amenability} 
	(extreme amenability in the domain
	of SUC flows). 
	We can similarly define \textbf{SUC-amenability} as the property of having a fixed point in 
	every affine SUC $G$-flow $Q$. Note that, the existence of a fixed point in $Q$ is equivalent to the existence of an  equicontinuous compact $G$-subspace of $Q$ (see \cite[Prop. 2.1]{GM-fixp12}). For SIN topological groups SUC-amenability and amenability are equivalent (see Remark \ref{r:abelian}).   	
\end{remark}


\begin{prop} \label{t:SUC-amenable} 
Let $G$ be a Polish group with metrizable $M(G)$. Then $G$ is SUC-amenable. 	
\end{prop} 
\begin{proof}
Let $Q$ be a compact affine	SUC $G$-system. 
It is enough to show that there exists an equicontinuous compact (minimal) $G$-subspace $X$ of $Q$. 
Choose any minimal $G$-subsystem $X$ of $Q$.
Then $X$ is SUC. Since $M(G)$ is metrizable and $X$ is a $G$-factor of $M(G)$, using again \cite[Prop. 14.1]{AKL}, there exists a point $x_0 \in X$ such that the orbit map $G \to Gx_0$ is open. Theorem \ref{t:stronger} guarantees that $X$ is equicontinuous.   
\end{proof}

\begin{prop}  \label{p:Ibarl}  
Every Roelcke precompact Polish group $G$ is SUC-amenable. 	
\end{prop}
\begin{proof}
By a result of Ibarlucia \cite[Theorem 2.9]{Ibar} 
for every Roelcke precompact Polish group $G$ holds  $\SUC(G)=\WAP(G)$. 
Choose any minimal $G$-subsystem $X$ of $Q$ and $x_0 \in X$. Then $X$ can be treated as a $G$-compactification $\nu \colon G \to X, g \mapsto gx_0$. 
Since $X$ is SUC, 
the corresponding algebra $\mathcal A_{\nu}$ of this compactification is a subalgebra of $\SUC(G)$. By our assumption,   
$\SUC(G)=\WAP(G)$. Hence, $X$ is WAP. Being minimal and WAP 
it is necessarily equicontinuous \cite[Cor. 6.11]{Menz}.    	
\end{proof}

 
 \begin{remark} \label{r:2subcl} 
 Propositions \ref{t:SUC-amenable} and \ref{p:Ibarl} provide two 
 sufficient conditions of SUC-amenability involving two important subclasses of Polish groups. These two classes are incomparable. 
 Indeed, by Pestov's result \cite{Pe-book}, the isometry group 
 $G:=\Iso(\U)$ of the Urysohn space $\U$ is extremely amenable (hence, $M(G)$ is metrizable) but $G$ is not Roelcke precompact. 
 On the other hand (answering a question from 
 \cite{MTT}), 
 there exist Roelcke precompact 
 Polish groups $G$ such that $M(G)$ is not metrizable (see, 
 \cite{EHN} 
 and 
 \cite{Kwiat}).  	
 \end{remark}
 
 \begin{remark} \label{r:EXAMPLE} 	
There exist nonamenable Polish groups (with metrizable $M(G)$) 
which is SUC-amenable but not SUC-extremely amenable. Indeed, 
take for example, the product $G:=H_+(\mathbb{T}) \times K$, where $K$ is a compact metrizable nontrivial group. Note that $M(G_1 \times G_2)=M(G_1) \times M(G_2)$ for Polish groups $G_1, G_2$ with metrizable $M(G_1), M(G_2)$ (see \cite[Example 3.5]{BZ}). Hence, $M(G)=M(H_+(\mathbb{T}))  \times K=\mathbb{T} \times K$. Now, observe that $G$ is SUC-amenable (Proposition \ref{t:SUC-amenable}), nonamenable (use
the fact that $H_+(\mathbb{T})$ is non-amenable), not SUC-extremely amenable ($M(G)$ has a nontrivial SUC $G$-factor $K$ without fixed points). 
 \end{remark}






\bibliographystyle{amsplain}

\end{document}